# On the elliptic sinh-Gordon equation with integrable boundary conditions II - Baker–Akhiezer theory.

10th December 2024

Graham Andrew Smith[1]

**Abstract:** We study finite-type solutions of the elliptic sinh-Gordon equation along the strip $\mathbb{R} \times [0, L]$ with Durham conditions on each boundary component. We determine necessary rationality criteria for these conditions to be satisfied on both components. These rationality criteria are analogous to the necessary and sufficient criteria for finite-type solutions to be periodic. However, in contrast to the periodic case, they are not in themselves sufficient for Durham conditions to be satisfied on both boundary components: certain additional criteria must also be satisfied at 4 specific points of the spectral curve. Nevertheless, in the special case where the Durham conditions imposed on the two boundary components are complementary to one another in a sense that we will clarify, these rationality criteria are, indeed, sufficient. This study is a necessary preliminary for any general construction of finite-type solutions of the elliptic sinh-Gordon equation subject to Durham boundary conditions.

**Classification AMS:** 53A10, 37K10

## 1 - Introduction.

**1.1 - Rationality conditions.** In this paper, we continue our study of real, singly-periodic solutions $u$ of the *elliptic sinh-Gordon equation*

$$u_{z\bar{z}} + \frac{1}{8}\sinh(2u) = 0 , \qquad (1.1)$$

defined over the strip $\mathbb{R} \times [0, L]$, and subject to the boundary conditions

$$u_y = Ae^u + Be^{-u} , \qquad (1.2)$$

where $A$ and $B$ are real and constant over each boundary component. Following [11], we will call (1.2) *Durham boundary conditions*. They are significant because they preserve the integrable nature of the problem, and were introduced independently by Corrigan in [9] and [10] and by Sklyanin in [23] and [24].

In [19], we showed that every periodic solution of (1.1)-(1.2) is of finite type, which means, grossly speaking, that it can be constructed via a certain ansatz out of polynomial data, which we call its *spectral data*. It is desirable for the purposes of construction and classification to express the boundary conditions (1.2) in terms of this spectral data. We address this partially in [19], and it is the aim of the current paper to complete this task.

We first recall how finite type solutions of the elliptic sinh-Gordon equation are constructed from spectral data (we refer the reader to [2] for a thorough treatment). We start with a genus $g$ hyperelliptic curve $\Sigma$, realized, as in Section 2.1, as a submanifold of $\hat{\mathbb{C}}^2$ with removable singularities at $0$ and $\infty$. We call $\Sigma$ the *spectral curve*, and we will suppose, furthermore, that it is invariant under the actions of the involutions

$$\sigma_1(\lambda, \nu) := (\lambda, -\nu) , \ \sigma_2(\lambda, \nu) := (\bar{\lambda}, \bar{\nu}) , \ \text{and } \sigma_3(\lambda, \nu) := \left(\frac{1}{\bar{\lambda}}, \frac{\bar{\nu}}{\bar{\lambda}^{g-1}}\right) .$$

In order to avoid technical complexity, we restrict ourselves in this paper to the case where $\Sigma$ has no branch points other than $0$ and $\infty$ over the real line $\mathbb{R}$ and the unit circle $\mathbb{S}^1$. In particular, this means that only spectral curves of even genus will be considered. The case of real branch points, and of spectral curves of odd genus, requires a more careful analysis, which we will discuss in later work.

---

[1] Departamento de Matemática, Pontifícia Universidade Católica do Rio de Janeiro (PUC-Rio), Rio de Janeiro, Brazil

Let $H^{1,0}(\Sigma)$ denote the space of first order holomorphic differentials of $\Sigma$, and recall that its dual space canonically identifies with $H^{0,1}(\Sigma)$. We denote by $\theta : H^{0,1}(\Sigma) \to \mathbb{C}$ its theta function and, given a point $z_0 \in H^{0,1}(\Sigma)$, we define the affine embedding $\phi[z_0] : \mathbb{R}^2 \to H^{0,1}(\Sigma)$ by

$$\phi[z_0](x,y) := z_0 + \mathrm{V}[(x+iy), -(x+iy)],$$

where the vector V, which depends $\mathbb{R}$-linearly on $(x+iy)$, is as in Section 2.2. It is a remarkable fact that, with $D \in H^{0,1}(\Sigma)$ as in Theorem 4.1 of [3], the function

$$u[z_0](x,y) := 2\mathrm{Ln}\left(\frac{\theta(\phi[z_0](x,y))}{\theta(\phi[z_0](x,y)+D)}\right), \tag{1.3}$$

is a solution of (1.1) and that, furthermore, all finite type, complex solutions of (1.1) are obtained in this manner from some spectral curve.

This construction is more clearly expressed in terms of the Jacobi variety of $\Sigma$. Recall that the integral homology group $H_1(\Sigma, \mathbb{Z})$ embeds as a cocompact lattice in $H^{0,1}(\Sigma)$, and that the quotient

$$\mathrm{Jac}(\Sigma) := H^{0,1}(\Sigma)/H_1(\Sigma, \mathbb{Z})$$

is a compact $g$-dimensional complex manifold, called the *Jacobi variety* of $\Sigma$. It turns out that if two points $z_0$ and $z_0'$ differ by an element of $H_1(\Sigma, \mathbb{Z})$, then the functions $u[z_0]$ and $u[z_0']$ coincide, so that solutions of (1.1) with spectral curve $\Sigma$ are parametrized by points of $\mathrm{Jac}(\Sigma)$. We henceforth identify points of $H^{0,1}(\Sigma)$ with their projections in $\mathrm{Jac}(\Sigma)$.

It is well-known that properties of the function $u[z_0]$ can be determined by properties of $z_0$. For example, concerning reality of solutions, it is well-known (see Section 3.3) that $\mathrm{Jac}(\Sigma)$ contains a compact, $g$-dimensional, totally real submanifold $\mathrm{Jac}_{\mathrm{re}}(\Sigma)$ such that $u[z_0]$ is real if and only if $z_0 \in \mathrm{Jac}_{\mathrm{re}}(\Sigma)$. We call this submanifold the *real locus* of $\Sigma$. It is the quotient of an affine subspace of $H^{0,1}(\Sigma)$, and is homeomorphic to $(\mathbb{S}^1)^g$. A non-trivial consequence is that every real, finite-type solution of (1.1) is smooth over the whole of $\mathbb{R}^2$ (see Remark 3.2).

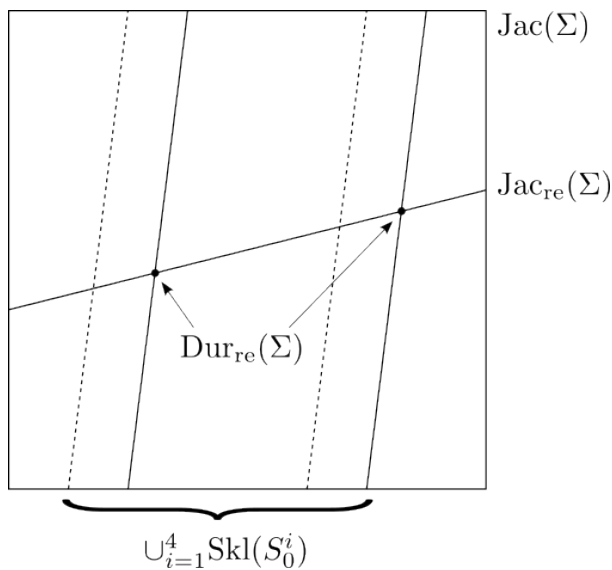

**Figure 1 - Components of the Jacobi variety** - The Jacobi variety $\mathrm{Jac}(\Sigma)$ is a quotient of $\mathbb{C}^g$ by a cocompact lattice. The *real locus* $\mathrm{Jac}_{\mathrm{re}}(\Sigma)$ is a real $g$-dimensional affine subspace, every point of which corresponds to a real solution of the sinh-Gordon equation. The four *Sklyanin subspaces* $\mathrm{Skl}(S_0^i)$ are complex $g/2$-dimensional affine subspaces. They intersect the real locus along four real $g/2$-dimensional affine subspaces which we call the *real Sklyanin subspaces*. The *real Durham locus*, here represented by two points, is a union of at most two real Sklyanin subspaces. Its points correspond to real solutions which satisfy the Durham boundary condition along $\mathbb{R} \times \{0\}$.

Our first main result concerns cohomological criteria for Durham boundary conditions to be satisfied over a single boundary component of $\mathbb{R} \times [0, L]$. This will require a number of technical definitions. We first define the *Sklyanin subset* $S$ of $\Sigma$ by

$$S := \{(\lambda, \nu) \in \Sigma \mid \mathrm{Det}(K(\lambda)) = 0\} \,, \tag{1.4}$$

where

$$K(\lambda) := \begin{pmatrix} 4A - 4B\lambda & \lambda - 1/\lambda \\ \lambda - 1/\lambda & 4A - 4B/\lambda \end{pmatrix} \tag{1.5}$$

denotes *Sklyanin's K-matrix*. We define a *special Sklyanin subset* of $\Sigma$ to be a subset $S_0$ of $S$ such that

$$S = S_0 \sqcup (\sigma_2 \sigma_3)(S_0) \,. \tag{1.6}$$

We say that two special Sklyanin subsets are *complementary* whenever their union is $S$, and we denote $S_0^* := S \setminus S_0$. When $|A| \neq |B|$, the Sklyanin subset has cardinality 8, each special Sklyanin subset has cardinality 4, there are precisely 4 such subsets, and any pairwise non-complementary family of such subsets has cardinality at most 2. The cases where $|A| = |B| \neq 0$ and $A = B = 0$ are discussed at the end of Section 4.1.

We define $\hat{\Psi} : H^{0,1}(\Sigma) \to H^{0,1}(\Sigma)$ by

$$\hat{\Psi}(z) := z - \sigma_2^* \sigma_3^* z \,. \tag{1.7}$$

Since $\hat{\Psi}$ sends $H_1(\Sigma, \mathbb{Z})$ to itself, it descends to a projection $\Psi$ from $\mathrm{Jac}(\Sigma)$ to itself. For each special Sklyanin subset $S_0$, we define the *real Sklyanin subspace* $\mathrm{Skl}_{\mathrm{re}}(S_0)$ by

$$\mathrm{Skl}_{\mathrm{re}}(S_0) := \{z \in \mathrm{Jac}_{\mathrm{re}}(\Sigma) \mid \Psi(z) = \mathcal{A}(S_0)\} \,, \tag{1.8}$$

where $\mathcal{A}$ denotes the Abel map normalized such that $\mathcal{A}(0) = 0$. There are trivially 4 such subspaces, each of which is a $(g/2)$-dimensional affine subspace of $\mathrm{Jac}_{\mathrm{re}}(\Sigma)$. It is straightforward to show that, for any special Sklyanin subset $S_0$,

$$\mathrm{Skl}_{\mathrm{re}}(S_0^*) = \mathrm{Skl}_{\mathrm{re}}(S \setminus S_0) = -\mathrm{Skl}_{\mathrm{re}}(S_0) \,. \tag{1.9}$$

We will say that two real Sklyanin subspaces are *complementary* whenever their Sklyanin subsets have this property. These objects are illustrated schematically in Figure 1.

**Theorem & Definition 1.1, The Durham locus**

There exists a subset $\mathrm{Dur}_{\mathrm{re}}(\Sigma)$ of $\mathrm{Jac}_{\mathrm{re}}(\Sigma)$, which is a union of pairwise non-complementary real Sklyanin subspaces, such that $u[z_0]$ is a real solution of (1.1) satisfying the Durham boundary condition (1.2) along $\mathbb{R} \times \{0\}$ if and only if $z_0 \in \mathrm{Dur}_{\mathrm{re}}(\Sigma)$. In particular, this set, which may be empty, has at most 2 connected components. We call $\mathrm{Dur}_{\mathrm{re}}(\Sigma)$ the *real Durham locus* of $\Sigma$.

**Remark 1.1.** Theorem 1.1 is proven in Section 4.3.

**Remark 1.2.** The most puzzling aspect of the Durham boundary conditions is the role played by the special Sklyanin subset $S_0$. Indeed, the Sklyanin subspace $\mathrm{Skl}_{\mathrm{re}}(S_0)$ is contained in the Durham locus if and only if a certain meromorphic function defined over $\Sigma$ assumes certain specific values over $S_0$ (see Lemma 4.4). We have not yet found a simple criterion for this condition to be satisfied. Nonetheless, and rather surprisingly, although this condition should apparently have codimension 4 in each Sklyanin subspace, it is in fact satisfied over the entire subspace as soon as it is satisfied at a single point.

Having homologically characterized the Durham boundary conditions, it is straightforward to determine a useful rationality condition for solutions of (1.1) to satisfy these criteria over both boundary components. For this, it suffices to observe, as we do in Lemma 4.14, that the *complementary real Durham locus*

$$\mathrm{Dur}_{\mathrm{re}}(\Sigma)^* := -\mathrm{Dur}_{\mathrm{re}}(\Sigma)$$

is none other than the real Durham locus of the *complementary boundary conditions*

$$u_y := -Ae^u - Be^u \tag{1.10}$$

over $\mathbb{R} \times \{0\}$. This allows us to address the case of greatest interest to us.

**Theorem 1.2, Rationality conditions**

If $z_0$ is a point of $\mathrm{Dur}_{\mathrm{re}}(\Sigma)$, and if

$$\Psi(2z_0) + \Psi(V[iL, -iL]) = 0 \ , \tag{1.11}$$

then the point

$$z_1 := z_0 + V[iL, -iL] \ , \tag{1.12}$$

lies on $\mathrm{Dur}_{\mathrm{re}}(\Sigma)^*$.

In particular, if $u := u[z_0]$ satisfies the Durham boundary conditions (1.2) along $\mathbb{R} \times \{0\}$, and if (1.11) holds, then $u$ satisfies the complementary Durham boundary conditions (1.10) along $\mathbb{R} \times \{L\}$.

**Remark 1.3.** Theorem 1.2 is proven in Section 4.3.

**Remark 1.4.** In [16], rationality conditions for periodic solutions are described in terms of the existence of a non-vanishing holomorphic function over $\Sigma$ with certain essential singularities at 0 and $\infty$. The rationality condition (1.11) readily admits an analogous reformulation.

**1.2 - Free boundary constant mean curvature surfaces.** Our interest in the sinh-Gordon equation is mainly motivated by its applications to the study of free boundary constant mean curvature (CMC) annuli in the unit ball in $\mathbb{R}^3$. Indeed, consider a conformal immersion $e : \mathbb{R} \times [0, \infty[ \to \mathbb{R}^3$ with constant Hopf differential $\phi := (1/4)dzdz$. It is a straightforward exercise of classical surface theory (see, for example, Section 1 of [19]) to show that $u$ has constant mean curvature equal to $(1/2)$ if and only if the conformal factor of its induced metric has the form $e^{2u}$, for some function $u : \mathbb{R} \times [0, \infty[ \to \mathbb{R}$ satisfying (1.1). Now let $\nu : \mathbb{R} \times \{0\} \to \mathbb{S}^2$ denote its outward-pointing, unit conormal vector field. An equally straightforward exercise shows that the boundary curve lies along a sphere of radius $R$, with $\nu$ making a constant angle $\theta \neq \pm\pi/2$ with the outward-pointing unit normal $n$ of that sphere if and only if (1.2) holds along $\mathbb{R} \times \{0\}$ with

$$A = \frac{\varepsilon}{\cos(\theta)R} + \frac{1}{2}\tan(\theta) \text{ and } B = \frac{1}{2}\tan(\theta) \ , \tag{1.13}$$

for some $\epsilon \in \{\pm 1\}$. For the reader's convenience, we prove this assertion in Appendix A.

The work [14] and [15] of Fraser–Schoen led to renewed interest in the study of free boundary minimal and CMC surfaces in the unit ball $B_1^3(0)$ of $\mathbb{R}^3$. Although it was initially expected for embedded free-boundary annuli in $B_1^3(0)$ to behave much like embedded tori in $\mathbb{S}^3$, a growing body of work now shows that these two classes differ in quite fundamental ways. The most striking of these results is the recent work [8] of Cerezo–Fernández–Mira, in which the authors construct embedded CMC annuli in the unit ball in $\mathbb{R}^3$ which are not surfaces of revolution. This solves in the negative an open problem posed by Wente in [26], and contrasts with the results [1] of Andrews–Li, [5] of Brendle, and [18] of Kilian–Schmidt for the case without boundary.

Cerezo–Fernández–Mira's work builds on the construction [25] of Wente. Thus, just as Wente's construction constitutes the starting point of the classification of CMC tori, initiated by Pinkall–Sterling in [22], and completed by Bobenko in [3], so too the work of Cerezo–Fernández–Mira emphasizes the interest in following a similar strategy to classify all immersed free boundary CMC annuli in the unit ball. The alignment conditions of Theorem 1.2 constitute a key component of this programme, and, in particular, a necessary step towards applying the construction techniques of [7], [12], and [17] to the free boundary case. We aim to address this in later work.

**1.3 - Structure of paper.** The paper is arranged as follows.

In Section 2, we introduce the spectral curve and we review the basic theories of abelian differentials, divisors, and potentials. In Section 3, we study Baker–Akhiezer theory. We review how this theory serves to construct complex solutions of (1.1), and we review the conditions required to ensure that these solutions are real.

The content of Sections 2 and 3, which constitutes about half of the paper, may be found scattered across the literature, and we have chosen to present it here in some detail, as we believe it will be of use to many readers.

Section 4 constitutes new material. Here we study how Durham boundary conditions are described in terms of affine conditions over the Jacobi variety. We introduce Sklyanin subspaces, which are complex $g/2$-dimensional affine subspaces of the Jacobi variety, as well as their real forms, which are real $g/2$-dimension affine subspaces of the real part of the Jacobi variety. We introduce the Durham locus, which consists of those points $z_0$ whose solution $u[z_0]$ satisfies (1.2) along $\mathbb{R} \times \{0\}$. In Section 4.1, we show that the Durham locus is contained in the union of all Sklyanin subspaces. In Section 4.2, we show that every connected component of the Durham locus is a Zariski open subset of one of the Sklyanin subspaces. In Section 4.3, we combine these properties to proof Theorems 1.1 and 1.2.

Finally, (1.13) is verified in Appendix A.

**1.4 - Acknowledgements.** The author is grateful to Martin Kilian for inspiring conversations, and also for many helpful comments made to earlier drafts of this paper.

## 2 - The spectral curve and its geometry.

**2.1 - The spectral curve.** We first review the algebraic framework used to construct finite-type solutions of the sinh-Gordon equation. Although most results of this section can be found scattered across the literature, we will provide proofs for the reader's convenience.

Let $g$ be an even integer, and let

$$\Delta(\lambda) := \sum_{i=0}^{2g} \lambda^i \Delta_i \ ,$$

be a polynomial of degree $2g$ having only simple roots, none of which lie on the real line $\mathbb{R}$, nor the unit circle $\mathbb{S}^1$. We will suppose that, for all $\lambda$,

$$\Delta(\bar{\lambda}) = \overline{\Delta(\lambda)} \text{ and } \Delta\left(\frac{1}{\bar{\lambda}}\right) = \frac{1}{\bar{\lambda}^{2g}} \overline{\Delta(\lambda)} \ .$$

We call these conditions the $\mathbb{R}$-*symmetry* and the $\mathbb{S}^1$-*symmetry* respectively. They are respectively equivalent to the conditions that, for all $i$,

$$\overline{\Delta}_i = \Delta_i \text{ and } \overline{\Delta}_{2g-i} = \Delta_i \ . \tag{2.1}$$

It will also be convenient to impose the normalization

$$\Delta_0 = \frac{1}{16} \ . \tag{2.2}$$

We will see presently that this imposes no restrictions on our construction. We will say that $\Delta$ is an *admissable spectral polynomial* whenever it satisfies the above conditions. In this case, its set Z of roots is symmetric under reflection in both the real line and the unit circle. In addition, the condition that no root lies on $\mathbb{R} \cup \mathbb{S}^1$ ensures that every orbit in Z of the group generated by these reflections has cardinality 4. Finally, since $\Delta$ has real coefficients and is non-vanishing over $\mathbb{R}$, the normalization (2.2) ensures that $\Delta$ is positive over this set.

We define the *finite part* of the spectral curve by

$$\Sigma^* := \left\{ (\lambda, \nu) \in \mathbb{C}^* \times \mathbb{C}^* \mid \nu^2 = -\frac{1}{\lambda}\Delta(\lambda) \right\} \ , \tag{2.3}$$

and we define the *spectral curve* $\Sigma$ to be its 2-point compactification. We verify that this curve has a single point at 0 and a single point at $\infty$, the restriction of $\nu$ to this curve has a simple pole at 0, and the restriction of $\nu\lambda^{1-g}$ has a simple pole at $\infty$. We denote by $\pi : \Sigma \to \hat{\mathbb{C}}$ the projection onto the first factor. We will only be concerned in the sequel with intrinsic properties of $\Sigma$. In particular, for $\alpha^2 \in \mathbb{R} \setminus \{0\}$, the biholomorphism $(\lambda, \nu) \mapsto (\lambda, \alpha\nu)$ preserves (2.1), which allows us to always recover the normalization (2.2), as indicated above.

The following simple decomposition result for meromorphic functions over $\Sigma$ will prove useful (c.f. Prop 1.10 of [21]).

**Lemma & Definition 2.1**

*Every meromorphic function $u : \Sigma \to \hat{\mathbb{C}}$ uniquely decomposes as*

$$u(\lambda, \nu) = v(\lambda) + \nu w(\lambda) , \qquad (2.4)$$

*where $v, w : \hat{\mathbb{C}} \to \hat{\mathbb{C}}$ are meromorphic functions which we respectively call its* even *and* odd *components.*

**Proof:** Indeed, we define $v, w : \hat{\mathbb{C}} \to \hat{\mathbb{C}}$ by

$$v(\lambda) := \frac{1}{2}(u(\lambda, \nu) + u(\lambda, -\nu)) \text{ and } w(\lambda) := \frac{1}{2\nu}(u(\lambda, \nu) - u(\lambda, -\nu)) .$$

It is straightforward to verify that these functions are well-defined and holomorphic, and the result follows. $\square$

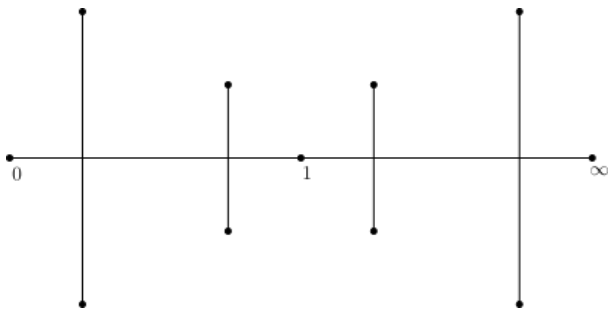

**Figure 2 - The branch cut** - The origin is joined to $\infty$ by the positive real axis. Each vertical line represents a circular arc, centred at 0 and passing through the positive real axis, joining a pair of conjugate roots of $\Delta$. This set is invariant under reflection in both $\mathbb{R}$ and $\mathbb{S}^1$.

The spectral curve is invariant under the action of the involutions

$$\sigma_1(\lambda, \nu) := (\lambda, -\nu) , \ \sigma_2(\lambda, \nu) := (\overline{\lambda}, \overline{\nu}) , \text{ and } \sigma_3(\lambda, \nu) := \left(\frac{1}{\overline{\lambda}}, \frac{\overline{\nu}}{\overline{\lambda}^{g-1}}\right) .$$

The first is known as the *hyperelliptic involution*, and we call the second and the third the $\mathbb{R}$-*involution* and the $\mathbb{S}^1$-*involution* respectively. They commute pairwise and generate an order 8 group isomorphic to $\mathbb{Z}_2^3$.

We now describe the explicit branch cut representation of $\Sigma$ that we will use. Let $B$ denote the branch cut in $\mathbb{C}$, symmetric under reflection in both $\mathbb{R}$ and $\mathbb{S}^1$, illustrated schematically in Figure 2. The complement of $\pi^{-1}(B)$ in $\Sigma$ consists of 2 simply-connected components, each of which projects biholomorphically onto $\mathbb{C} \setminus B$. These components are identified by the sign of $\nu$ over the negative real axis. Indeed, by (2.2), $\nu$ is always real and non-zero over this axis. We call the component of $\Sigma \setminus \pi^{-1}(B)$ where $\nu$ is positive over this axis the *upper sheet* and we call the other component the *lower sheet*. We represent curves in the upper and lower sheets respectively by solid and broken lines.

**Lemma 2.2**

*(1) The hyperelliptic involution $\sigma_1$ exchanges sheets;*

*(2) the $\mathbb{R}$-involution $\sigma_2$ preserves sheets; and*

*(3) the $\mathbb{S}^1$-involution $\sigma_3$ exchanges sheets.*

**Proof:** It suffices to determine the actions of these involutions at $\pi^{-1}(\{-1\})$, bearing in mind that the genus $g$ is even. $\square$

We now describe the system of canonical generators of $H_1(\Sigma)$ that we will use. Recall that this consists of a pair of ordered sets $(\alpha_1, \cdots, \alpha_g)$ and $(\beta_1, \cdots, \beta_g)$ of oriented, simple closed curves in $\Sigma$, respectively called $\alpha$-*cycles* and $\beta$-*cycles*, such that each $\alpha$-cycle crosses its corresponding $\beta$-cycle exactly once from left to right, and no other pair of these curves intersects at any other point. We illustrate the first two elements of the sets of $\alpha$-cycles and $\beta$-cycles that we will use in Figure 3, leaving the reader to infer how the remaining elements are defined.

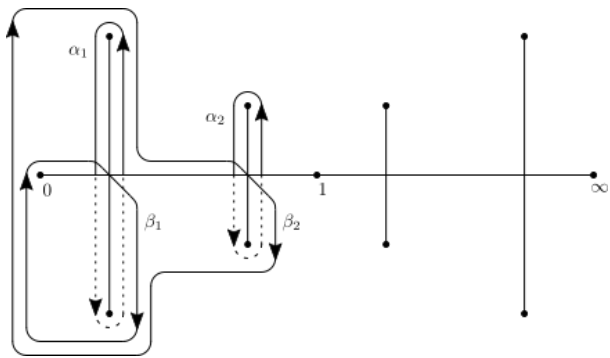

**Figure 3 - The $\alpha$- and $\beta$-cycles** - Each $\alpha$-cycle turns anticlockwise around a vertical branch cut, spending equal time in each sheet. The $\beta$-cycles all lie in the upper sheet, and are nested in order to avoid intersections.

It is straightforward to determine the actions of the involutions $\sigma_1$, $\sigma_2$, and $\sigma_3$ on the homology classes of the canonical generators. Indeed, their actions on the $\alpha_i$'s are

$$\sigma_1 \cdot \alpha_i = -\alpha_i \ , \ \sigma_2 \cdot \alpha_i = \alpha_i \ , \text{ and } \sigma_3 \cdot \alpha_i = \alpha_{(g+1)-i} \ , \tag{2.5}$$

whilst their actions on the $\beta_i$'s are

$$\sigma_1 \cdot \beta_i = -\beta_i \ , \ \sigma_2 \cdot \beta_i = -\beta_i + \alpha_i \ , \text{ and } \sigma_3 \cdot \beta_i = -\beta_{(g+1)-i} + \alpha_{(g+1)-i} \ . \tag{2.6}$$

**2.2 - Abelian differentials.** We now review the effect of the actions of $\sigma_1$, $\sigma_2$, and $\sigma_3$ on abelian differentials and related objects over $\Sigma$. We first address abelian differentials of the first kind, that is, holomorphic differentials. Let $\zeta_1, \cdots, \zeta_g$ be a basis of $H^{1,0}(\Sigma)$, normalized in the sense that, for all $(i,j)$,

$$\int_{\alpha_i} \zeta_j = \delta_{ij} \ ,$$

where $\delta_{ij}$ denotes the Kronecker delta function.

**Lemma 2.3**

*For all $i$,*

$$\sigma_1^* \zeta_i = -\zeta_i \ , \ \sigma_2^* \zeta_i = \overline{\zeta}_i \ , \text{ and } \sigma_3^* \zeta_i = \overline{\zeta}_{(g+1)-i} \ . \tag{2.7}$$

**Proof:** We only prove the second relation, as the proofs of the remaining two are similar. Since $\sigma_2$ is antiholomorphic, $\sigma_2^* \overline{\zeta}_i$ is a holomorphic differential. By (2.5), for all $j$,

$$\int_{\alpha_j} \sigma_2^* \overline{\zeta}_i = \int_{\sigma_2 \cdot \alpha_j} \overline{\zeta}_i = \int_{\alpha_j} \overline{\zeta}_i = \delta_{ji} \ ,$$

so that $\sigma_2^* \overline{\zeta}_i = \zeta_i$. The result now follows upon conjugating this identity. $\square$

The *period matrix* of $\Sigma$ is defined by

$$\Pi_{ij} := \int_{\beta_i} \zeta_j \ . \tag{2.8}$$

Recall (see [13]) that this matrix is symmetric with positive-definite imaginary part.

**Lemma 2.4**

For all $(i,j)$,
$$\operatorname{Re}(\Pi_{ij}) = \frac{1}{2}\delta_{ij} \ . \tag{2.9}$$

**Proof:** Indeed, for all $(i,j)$, by (2.6) and (2.7),
$$\overline{\Pi}_{ij} = \int_{\beta_i} \overline{\zeta}_j = \int_{\beta_i} \sigma_2^*\zeta_j = \int_{\sigma_2\cdot\beta_i} \zeta_j = -\int_{\beta_i} \zeta_j + \int_{\alpha_i} \zeta_j = -\Pi_{ij} + \delta_{ij} \ ,$$
and the result follows. □

**Lemma 2.5**

For all $(i,j)$,
$$\Pi_{i,(g+1)-j} = \Pi_{(g+1)-i,j} \ . \tag{2.10}$$

**Proof:** Indeed, for all $(i,j)$, by (2.5) and (2.7),
$$\overline{\Pi}_{i,(g+1)-j} = \int_{\beta_i} \overline{\zeta}_{(g+1)-j} = \int_{\beta_i} \sigma_3^*\zeta_j = \int_{\sigma_3\cdot\beta_i} \zeta_j = -\int_{\beta_{(g+1)-i}} \zeta_j + \int_{\alpha_{(g+1)-i}} \zeta_j = -\Pi_{(g+1)-i,j} + \delta_{(g+1)-i,j} \ ,$$
and the result now follows by (2.9). □

Recall that the *Jacobi variety* of $\Sigma$ is defined by
$$\operatorname{Jac}(\Sigma) := \mathbb{C}^g/\Lambda \ ,$$
where $\Lambda$ is the lattice generated by the canonical basis $e_1,\cdots,e_g$, together with the vectors $f_1,\cdots,f_g$ defined by
$$f_i := \sum_{j=1}^{g} \Pi_{ij}e_j \ .$$

Let $\mathcal{A}:\Sigma\to\operatorname{Jac}(\Sigma)$ denote the *Abel map* of $\Sigma$ based at $0$, that is
$$\mathcal{A}(p)_i := \int_0^p \zeta_i \ . \tag{2.11}$$

**Lemma 2.6**

For all $i$,
$$(\mathcal{A}\circ\sigma_1)_i = -\mathcal{A}_i \ , \ (\mathcal{A}\circ\sigma_2)_i = \overline{\mathcal{A}}_i \ , \text{ and } (\mathcal{A}\circ\sigma_3)_i = \mathcal{A}(\infty)_i + \overline{\mathcal{A}}_{(g+1)-i} \ . \tag{2.12}$$

**Proof:** It suffices to prove the second identity, as the proofs of the first and third are similar. Since $\sigma_2(0)=0$, for all $p\in\Sigma$,
$$(A\circ\sigma_2)(p)_i = \int_0^{\sigma_2(p)} \zeta_i = \int_{\sigma_2(0)}^{\sigma_2(p)} \zeta_i \ .$$

Thus, by (2.7),
$$(A\circ\sigma_2)(p)_i = \int_{\sigma_2(0)}^{\sigma_2(p)} \sigma_2^*\overline{\zeta}_i = \int_0^p \overline{\zeta}_i = \overline{\mathcal{A}(p)}_i \ ,$$
as desired. □

We now address abelian differentials of the second kind, that is, meromorphic differentials, every one of whose poles has vanishing residue. We say that an abelian differential of the second kind $\xi$ is *normalized* whenever, for all $i$,
$$\int_{\alpha_i} \xi = 0 \ .$$
Since any normalized holomorphic differential vanishes, it follows that $\xi$ is uniquely determined by this normalization and the principal parts of its poles.

Let $P$ and $Q$ be polynomials. We define $\Omega := \Omega[P,Q]$ to be the unique normalised abelian differential of the second kind with poles at $0$ and $\infty$ given respectively by
$$\Omega(\lambda,\nu) = d\bigl(\nu P(1/\lambda)\bigr) + \mathrm{O}(1) \text{ and } (\sigma_3^*\Omega)(\lambda,\nu) = d\bigl(\overline{\nu Q}(1/\overline{\lambda})\bigr) + \mathrm{O}(1) \ .$$
We define the vector $\mathrm{V}[P,Q] \in \mathbb{C}^g$ such that, for all $i$,
$$\mathrm{V}[P,Q]_i := \int_{\beta_i} \Omega[P,Q] \ . \tag{2.13}$$
Trivially, both $\Omega$ and $\mathrm{V}$ depend $\mathbb{R}$-linearly on $(P,Q)$, $\mathbb{C}$-linearly on $P$, and $\mathbb{C}$-antilinearly on $Q$. The following result will allow us later to determine which differentials yield non-trivial flows.

**Lemma 2.7**

*If $P$ has degree no greater than $(g-1)$, and if $V[P,0]$ vanishes, then so too does $P$.*

**Remark 2.1.** In particular, the map $P \mapsto \mathrm{V}[P,0]$ defines a $\mathbb{C}$-linear isomorphism from the space of complex degree $(g-1)$ polynomials into $\mathbb{C}^g$.

**Proof:** Indeed, when $\mathrm{V}[P,0] = 0$, the normalized abelian differential $\Omega[P,0]$ integrates to a meromorphic function $u : \Sigma \to \hat{\mathbb{C}}$ having a pole only at $0$ which, furthermore, has principal part $\nu P(1/\lambda)$. We claim that $u$ is constant. Indeed, let $v,w : \hat{\mathbb{C}} \to \hat{\mathbb{C}}$ denote respectively the even and odd components of $u$ as in (2.4), so that
$$u = v + \nu w \ .$$
Since the principal part of the pole of $u$ at $0$ only has terms of odd order in $\nu$, it follows that $v$ has no poles, and is therefore constant. We now claim that $w$ vanishes. Indeed, suppose the contrary. Since $\Sigma$ only has simple branch points, $w$ has no pole at any point of $\mathbb{C}^*$. Furthermore, near $\infty$,
$$|\nu| = \mathrm{O}\bigl(|\lambda|^{g-1/2}\bigr) \text{ and } u = \mathrm{O}(1) \ ,$$
so that $w$ has a zero of order at least $g$ at this point, and therefore also a pole of order at least $g$ at $0$. However, since the principal part of $w$ at $0$ is equal to $P(1/\lambda)$, the order of this pole is at most $(g-1)$, which is absurd. It follows that $w$ vanishes as asserted, and so too does $\Omega = du = dv$, and therefore also $P$. This completes the proof. $\square$

**Lemma 2.8**

*For all $P,Q$,*
$$\Omega[\overline{P},\overline{Q}] = \sigma_2^*\overline{\Omega[P,Q]} \text{ and } \Omega[Q,P] = \sigma_3^*\overline{\Omega[P,Q]} \ . \tag{2.14}$$

**Proof:** We will only prove the first relation, as the proof of the second is identical. Denote $\Omega := \Omega[P,Q]$. Since $\sigma_2$ is anti-holomorphic, $\sigma_2^*\overline{\Omega}$ is an abelian differential of the second kind with poles at $0$ and $\infty$ respectively given by
$$(\sigma_2^*\overline{\Omega})(\lambda,\nu) = d\bigl(\nu \overline{P}(1/\lambda)\bigr) + \mathrm{O}(1) \text{ and } (\sigma_3^*\sigma_2^*\overline{\Omega})(\lambda,\nu) = d\bigl(\overline{\nu}Q(1/\overline{\lambda})\bigr) + \mathrm{O}(1) \ .$$
Furthermore, bearing in mind (2.5), for all $i$,
$$\int_{\alpha_i} \sigma_2^*\overline{\Omega} = \int_{\sigma_2\cdot\alpha_i} \overline{\Omega} = \int_{\alpha_i} \overline{\Omega} = 0 \ ,$$
so that $\sigma_2^*\overline{\Omega}$ is normalized, and the result follows. $\square$

**Lemma 2.9**

For all $P, Q$, and for all $i$,

$$V[\overline{P}, \overline{Q}]_i = -\overline{V[P,Q]}_i \text{ and } V[Q, P]_i = \overline{V[P,Q]}_{(g+1)-i} . \tag{2.15}$$

**Proof:** We will only prove the first relation, as the proof of the second is identical. By (2.14), for all $i$,

$$V[\overline{P}, \overline{Q}]_i = \int_{\beta_i} \Omega[\overline{P}, \overline{Q}] = \int_{\beta_i} \sigma_2^* \overline{\Omega[P,Q]} = \int_{\sigma_2 \cdot \beta_i} \overline{\Omega[P,Q]} .$$

Thus, by (2.5),

$$V[\overline{P}, \overline{Q}]_i = -\int_{\beta_i} \overline{\Omega[P,Q]} + \int_{\alpha_i} \overline{\Omega[P,Q]} = -\overline{V[P,Q]}_i ,$$

as desired. $\square$

**2.3 - Potentials and their divisors.** We now introduce potentials of the spectral curve. These are simple Laurent polynomials[2], taking values in the Lie algebra $\mathfrak{sl}(2,\mathbb{C})$, and satisfying certain properties that will be described presently. Potentials play two roles in this theory: they are used to characterize non-special divisors supported over $\Sigma^*$, but, above all, they serve to construct complex solutions of (1.1), as we will see in Section 3.

We first review the properties of divisors of meromorphic functions of $\Sigma$. Given such a meromorphic function $\phi : \Sigma \to \hat{\mathbb{C}}$ and a point $z \in \hat{\mathbb{C}}$, we denote by $[\phi]_z$ the pole divisor of $\phi$ when $z = \infty$, and the zero divisor of $(\phi - z)$ otherwise. For any divisor $D$ over $\Sigma$, we call its restriction $D_{\text{fin}}$ to $\Sigma^*$ its *finite part*. In particular, we denote the finite part of $[\phi]_z$ by $[\phi]_{z,\text{fin}}$.

**Lemma 2.10**

Let $D$ be a positive, degree $g$ divisor supported over $\Sigma^*$. If $\phi : \Sigma \to \hat{\mathbb{C}}$ is a meromorphic function such that

$$[\phi]_\infty \leqslant D , \tag{2.16}$$

then $\phi$ depends only on $\lambda$.

**Proof:** Let $u$ and $v$ denote respectively the even and odd components of $\phi$ as in (2.4), so that

$$\phi(\lambda, \nu) = u(\lambda) + \nu v(\lambda) .$$

We show that $v$ vanishes. Indeed, suppose first that $v$ is non-constant, and recall that

$$v(\lambda) = \frac{1}{2\nu} \big(\phi(\lambda, \nu) - \phi(\lambda, -\nu)\big) .$$

Since $\Sigma$ only has simple branch points, the pole divisor of $v$ is bounded above by $\pi(D)$. In particular, $v$ has at most $g$ poles over $\hat{\mathbb{C}}$. However, since $\nu$ has a pole at 0, and since $\phi = O(1)$ near this point, $v$ has a root at 0. Likewise, since, near $\infty$,

$$|\nu| = O(|\lambda|^{g-1/2}) \text{ and } \phi = O(1) ,$$

the function $v$ has a root of order at least $g$ at this point. The function $v$ therefore has at least $(g+1)$ roots over $\hat{\mathbb{C}}$, which is absurd, since it has as many roots as poles. It follows that $v$ is constant, and since it vanishes at $\infty$, it vanishes everywhere, as desired. $\square$

This allows us to determine whenever a degree $g$ divisor supported over $\Sigma^*$ is non-special.

---

[2] We will say that a Laurent polynomial is *simple* whenever its lowest order term has order $(-1)$.

**Lemma 2.11**

Let $D$ be a positive, degree $g$ divisor supported over $\Sigma^*$. $D$ is special if and only if there exists a divisor $D'$ over $\hat{\mathbb{C}}$ such that

$$\pi^* D' \leqslant D . \tag{2.17}$$

**Remark 2.2.** That is, $D$ is non-special whenever the preimage of any point of $\mathbb{C}^*$ meets $D$ at a single point and, moreover, any point of $D$ that is also a branch point of $\Sigma$ has multiplicity at most 1.

**Proof:** Recall that $D$ is special if and only if there exists a non-constant meromorphic function $\phi : \Sigma \to \hat{\mathbb{C}}$ such that $[\phi]_\infty \leqslant D$. The result now follows upon taking $\phi(\lambda, \nu) = u(\lambda)$, where $u$ is a rational function with pole divisor $D'$. $\square$

Let $T^+, T^- \subseteq \mathfrak{sl}(2, \mathbb{C})$ denote respectively the Lie subalgebras of strictly upper- and lower-triangular matrices. We define a *potential* for $\Sigma$ to be a polynomial

$$\xi := \sum_{k=-1}^{g} \lambda^k \xi_k := \begin{pmatrix} \omega & \sigma \\ \tau & -\omega \end{pmatrix} := \sum_{k=-1}^{g} \lambda^k \begin{pmatrix} \omega_k & \sigma_k \\ \tau_k & -\omega_k \end{pmatrix} \tag{2.18}$$

such that

$$\xi_{-1} \in T^+ , \; \xi_g \in T^- , \; \text{and} \; \mathrm{Det}(\xi(\lambda)) = \frac{1}{\lambda}\Delta(\lambda) . \tag{2.19}$$

It will also be convenient to impose the normalization

$$\sigma_{-1} \sigma_{g-1} = \frac{1}{16} . \tag{2.20}$$

By (2.19) and the definition of $\Sigma$, for all $\lambda$, the eigenvalues of $\xi(\lambda)$ are $\pm \nu$. Note also that

$$\omega^2 + \sigma\tau = -\mathrm{Det}(\xi(\lambda)) = -\frac{1}{\lambda}\Delta(\lambda) = \nu^2 , \tag{2.21}$$

which yields the useful relation

$$(\nu - \omega)(\nu + \omega) = \sigma\tau . \tag{2.22}$$

We now study how potentials of $\Sigma$ are related to non-special degree $g$ divisors supported over $\Sigma^*$. Let $\xi$ be a potential of $\Sigma$. We define its *eigenline bundle* to be the complex subbundle of $\Sigma \times \mathbb{C}^2$ whose fibre at every non-branch point $(\lambda, \nu)$ is the eigenline of $\xi(\lambda)$ with eigenvalue $\nu$. We identify $\mathbb{CP}(1)$ with $\hat{\mathbb{C}}$ by identifying the point $[z : w]$ with the point $w/z$. The image of the eigenline bundle under this identification is a meromorphic function $\phi : \Sigma \to \hat{\mathbb{C}}$ which we call the *projectivized eigenline* of $\xi$.

**Lemma 2.12**

The projectivized eigenline of $\xi$ is given by

$$\phi(\lambda, \nu) = \frac{\nu - \omega(\lambda)}{\sigma(\lambda)} . \tag{2.23}$$

**Proof:** Indeed, by (2.18), for every non-branch point $(\lambda, \nu)$,

$$\xi(\lambda) \begin{pmatrix} 1 \\ \phi(\lambda, \nu) \end{pmatrix} = \begin{pmatrix} \nu \\ \nu\phi(\lambda, \nu) \end{pmatrix} ,$$

and the result follows. $\square$

**Lemma & Definition 2.13**

*The projectivized eigenline $\phi$ has a simple zero at $0$, a simple pole at $\infty$, and the finite part $D := [\phi]_{\infty,\text{fin}}$ of its pole divisor is a degree $g$ non-special divisor. We call $D$ the divisor of $\xi$.*

**Proof:** Since $\xi(0)$ is upper triangular, and since $\Delta(0) \neq 0$, $\omega_{-1} = 0$ and $\sigma_{-1} \neq 0$. It follows that, near $0$,

$$\phi(\lambda, \nu) = \frac{\lambda \nu}{\sigma_1} + \mathrm{O}(\lambda) \ .$$

so that $\phi$ has a simple root at $0$. In a similar manner, we show that, near $0$,

$$\phi\left(\frac{1}{\lambda}, \frac{\nu}{\lambda^{g-1}}\right) = \frac{\nu}{\sigma_{g-1}} + \mathrm{O}(1) \ ,$$

so that $\phi$ has a simple pole at $\infty$.

We now study the poles of $\phi$ over $\Sigma^*$. These poles are a subset of the set of roots of $\sigma$ over $\Sigma^*$. Let $\lambda_0 \in \mathbb{C}^*$ be a root of $\sigma$. Suppose first that $\lambda_0$ is not a root of $\Delta$. Let $\nu_0 \neq 0$ be such that $(\lambda_0, \pm \nu_0)$ are the points of $\Sigma$ lying above $\lambda_0$. By (2.22), we may assume that $\nu_0 + \omega(\lambda_0) = 0$ and, since $\nu_0 - \omega_0(\lambda_0) \neq 0$, the order of the root of $\nu + \omega(\lambda)$ is at least equal to that of $\sigma$ at $\lambda_0$. It then follows that $\phi$ has no pole at $(\lambda_0, -\nu_0)$ and a pole at $(\lambda_0, \nu_0)$ of order equal to that of the root of $\sigma$ at $\lambda_0$. Suppose now that $\lambda_0$ is a root of $\Delta$. Then, by (2.22), $\omega(\lambda_0) = 0$ and, since the root of $\Delta$ at $\lambda_0$ is simple, near this point,

$$\sigma\tau = \nu^2 + \mathrm{O}(\nu^3) = c(\lambda - \lambda_0) + \mathrm{O}(\nu^3) \ ,$$

for some $c \neq 0$. Consequently, $\sigma\tau$, and thus a fortiori $\sigma$, has at most a simple root at $\lambda_0$. It then follows by (2.23) that $\phi$ has a simple pole at this point. Since this accounts for all the poles of $\phi$ over $\Sigma^*$, we conclude that

$$\mathrm{Deg}([\phi]_{\infty,\text{fin}}) = \mathrm{Deg}(\lambda\sigma) = g \ .$$

Finally, by Lemma 2.11, this divisor is non-special, and this completes the proof. $\square$

The converse of this result also holds.

**Lemma 2.14**

*Every degree $g$ non-special divisor supported over $\Sigma^*$ is the divisor of some potential.*

**Remark 2.3.** It also follows from the proof of Lemma 2.14 that the operation sending potentials to their divisors is a $2 - 1$ map. We will see in Section 3.3 how this ambiguity is resolved for real solutions.

**Proof:** Let $D$ be a degree $g$ non-special divisor supported over $\Sigma^*$. Let $D^*$ denote the image of $D$ under the hyperelliptic involution. Let $\sigma$ be a simple Laurent polynomial of degree $(g-1)$ vanishing at $\pi(D)$, and satisfying (2.20). Observe that $\sigma$ is uniquely defined by this condition up to a change of sign. Let $\omega$ be the unique polynomial of degree $(g-1)$ such that, at every point $(\lambda, \nu)$ of $D^*$,

$$\omega(\lambda) = \nu \ .$$

We verify that $\phi := (\nu - \omega)/\sigma$ has a simple zero at $0$, a simple pole at $\infty$, and the finite part of its pole divisor is $D$. Furthermore, the function $(\nu - \omega)(\nu + \omega)$ is invariant under the hyperelliptic involution, and is thus equal to $p(\lambda)$, for some simple Laurent polynomial $p$. Since the roots of $\sigma$ are also roots of $p$, there exists a unique polynomial $\tau$ such that

$$\sigma\tau = p = (\nu - \omega)(\nu + \omega) \ .$$

We readily verify that $\sigma$, $\tau$ and $\omega$ are the components of a potential $\xi$ with divisor $D$, and this completes the proof. $\square$

Lemmas 2.13 and 2.14 justify the following definition.

**Definition 2.15**

*We say that a point $z \in \mathrm{Jac}(\Sigma)$ is a* potential point *whenever it is the image under the Abel map of a non-special degree $g$ divisor supported over $\Sigma^*$.*

We conclude this section by reviewing certain symmetries of the space of potentials. Let $\xi$ be a potential with divisor $D$. The function $\xi' := -\xi$ is trivially also a potential with divisor $D' := \sigma_1(D)$ and projectivized eigenline

$$\phi'(\lambda, \nu) := \frac{-\nu - \omega}{\sigma} \ .$$

Likewise, we readily verify that the gauge-transformed function

$$\xi'' := \begin{pmatrix} -\omega & \sigma \\ \tau & \omega \end{pmatrix} = \left[ \begin{pmatrix} 0 & 1 \\ 1 & 0 \end{pmatrix} \xi \begin{pmatrix} 0 & 1 \\ 1 & 0 \end{pmatrix} \right]^t, \tag{2.24}$$

is also a potential, with divisor $D'' := \sigma_1(D)$, and projectivized eigenline

$$\phi'(\lambda, \nu) := \frac{\nu + \omega}{\sigma} \ .$$

## 3 - Baker–Akhiezer theory.

**3.1 - Baker–Akhiezer flows.** Recall that integrable systems are in general characterized by possessing countable families of Poisson-commuting invariants of motion which, by Noether's theorem, identify with countable families of commuting flows. In this paper, we will be concerned with realizations of these flows in three different spaces, namely in the space of potentials, in the space of non-special degree $d$ divisors supported on $\Sigma^*$, and in the Jacobi variety. The lexicon allowing us to pass between any two is provided by Baker–Akhiezer theory. We present here an overview of the main steps, referring the reader to [2] and [20] for a thorough treatment of these ideas.

We first introduce Baker–Akhiezer flows. These are flows in the space of non-special degree $g$ divisors which correspond to straight-line flows in the Jacobi variety. They are defined as follows. Let $D$ be a non-special, degree $g$ divisor supported over $\Sigma^*$, and let $P$ and $Q$ be polynomials. Following [2], for sufficiently small $t$, there exist meromorphic functions $f, g : \Sigma^* \to \hat{\mathbb{C}}$, unique up to sign, both with pole divisor $D$, and with essential singularities at $0$ and $\infty$ given by

$$\begin{aligned}
f_t(\lambda, \nu) &= \mathrm{Exp}\Big( t\nu P\Big(\frac{1}{\lambda}\Big)\Big)(\rho_t + \mathrm{O}(\lambda\nu)) \ , \\
(f_t \circ \sigma_3)(\lambda, \nu) &= \mathrm{Exp}\Big( t\overline{\nu}\overline{Q}\Big(\frac{1}{\overline{\lambda}}\Big)\Big)\Big(\frac{1}{\rho_t} + \mathrm{O}(\lambda\nu)\Big) \ , \\
g_t(\lambda, \nu) &= \mathrm{Exp}\Big( t\nu P\Big(\frac{1}{\lambda}\Big)\Big)(4i\lambda\nu\sigma_t + \mathrm{O}(\lambda)) \ , \text{ and} \\
(g_t \circ \sigma_3)(\lambda, \nu) &= \mathrm{Exp}\Big( t\overline{\nu}\overline{Q}\Big(\frac{1}{\overline{\lambda}}\Big)\Big)\Big(\frac{i}{4\overline{\lambda}\overline{\nu}\sigma_t} + \mathrm{O}(1)\Big) \ ,
\end{aligned} \tag{3.1}$$

where $\rho_t$ and $\sigma_t$ are non-vanishing functions, depending only on $t$. The functions $f$ and $g$ are known as *Baker–Akhiezer functions*.

Denote

$$\phi_t := \frac{g_t}{f_t} \ .$$

Since the exponential factors of the essential singularities of $f$ and $g$ cancel, it is straightforward to show that $\phi_t$ has a simple root at $0$ and a simple pole at $\infty$. For all $t$, let $D_t$ denote the finite part of its pole divisor, as defined in Section 2.3. We call $(D_t)_{t \in ]-\delta, \delta[}$ the *Baker–Akhiezer flow passing through $D$ defined by $(P, Q)$*. By the construction of $f$ and $g$ given in [2], for all $t$,

$$\mathcal{A}(D_t) = \mathcal{A}(D) + t\mathrm{V}[P, Q] \ ,$$

where V[P,Q] is given by (2.13). The Baker–Akhiezer flow defined by $(P,Q)$ is thus sent by the Abel map to a straight-line flow in the Jacobi variety of constant velocity equal to V[P,Q].

We now study the corresponding flow in the space of potentials. Consider the matrix-valued function

$$\xi_t := \frac{2\nu}{\phi_t - (\sigma_1^*\phi_t)} \begin{pmatrix} -(\phi_t + \sigma_1^*\phi_t)/2 & 1 \\ -\phi_t(\sigma_1^*\phi_t) & (\phi_t + \sigma_1^*\phi_t)/2 \end{pmatrix} . \tag{3.2}$$

Since this function is invariant under the hyperelliptic involution, it only depends on $\lambda$. We will now show that it defines a flow in the space of potentials corresponding to $(D_t)_{t \in ]-\delta,\delta[}$.

**Lemma 3.1**

*For sufficiently small $t$, the matrix-valued function $\xi_t$ is a potential of $\Sigma$ with projectivized eigenline $\phi_t$ and divisor $D_t$.*

**Proof:** Indeed, let $u$ and $v$ denote respectively the odd and even components of $\phi_t$, so that

$$\phi_t = \nu u + v .$$

In particular,

$$u = \frac{1}{2\nu}(\phi_t - \sigma_1^*\phi_t) , \text{ and}$$
$$v = \frac{1}{2}(\phi_t + \sigma_1^*\phi_t) .$$

Denote also

$$w := \phi_t(\sigma_1^*\phi_t) ,$$

and observe that this function also only depends on $\lambda$. By (3.1), near 0,

$$\begin{aligned}\phi_t(\lambda,\nu) &= \frac{4i\lambda\nu\sigma_t}{\rho_t} + O(\lambda) , \text{ and} \\ (\phi_t \circ \sigma_3)(\lambda,\nu) &= \frac{i\rho_t}{4\overline{\lambda}\overline{\nu}\sigma_t} + O(1) .\end{aligned} \tag{3.3}$$

We now analyse the roots and poles of $u$, $v$ and $w$. Since $\Sigma$ has simple branch points, and since the finite part of the pole divisor of $\phi$ is non-special, $u$ has exactly $\mathrm{Deg}(D) = g$ poles over $\mathbb{C}^*$. Next, by (3.3), near 0,

$$u(\lambda) = \frac{4i\lambda\sigma_t}{\rho_t} + O(\lambda^{3/2}) , \text{ and}$$
$$u\left(\frac{1}{\overline{\lambda}}\right) = \frac{i\rho_t\overline{\lambda}^{g-1}}{4\overline{\lambda}\overline{\nu}^2\sigma_t} + O(\lambda^{g-1/2}) .$$

This function therefore has a simple root at 0, a root of order $(g-1)$ at $\infty$, and, since it has the same number of roots as poles, it has no other root over $\mathbb{C}$. It follows that $\sigma := 1/u$ is a simple Laurent polynomial of degree $(g-1)$ such that

$$\sigma_{-1} = \frac{\rho_t}{4i\sigma_t} \neq 0 , \tag{3.4}$$

and, by (2.2),

$$\sigma_{g-1} = \lim_{(\lambda,\nu) \to 0} \frac{4\overline{\lambda}\overline{\nu}^2\sigma_t}{i\rho_t} = \frac{-4\Delta_0\sigma_t}{i\rho_t} = \frac{-\sigma_t}{4i\rho_t} , \tag{3.5}$$

so that the normalization condition (2.20) holds.

Since the finite part of the pole divisor of $\phi$ is non-special, $v$ also has at most $\mathrm{Deg}(D) = g$ poles over $\mathbb{C}^*$, which coincide with poles of $u$, and so $v/u$ is regular over $\mathbb{C}^*$. Furthermore, by (3.3) again, near 0,

$$v(\lambda) = O(\lambda) , \text{ and}$$
$$v\left(\frac{1}{\overline{\lambda}}\right) = O(1) ,$$

so that $\omega := -v/u$ is a polynomial of degree $(g-1)$.

Finally, since the pole divisor of $\phi$ is non-special, $w$ also has exactly $\text{Deg}(D) = g$ poles over $\mathbb{C}^*$, which coincide with those of $u$. Furthermore, by (3.3) again, near $0$,

$$w(\lambda) = O(\lambda) , \text{ and}$$

$$w\left(\frac{1}{\bar\lambda}\right) = \frac{\rho_t^2}{16\bar\lambda^2 \bar\nu^2 \sigma_t^2} + O(1) ,$$

so that $\tau := -w/u$ is a polynomial of degree $g$.

Finally, we verify that, for all $(\lambda, \nu)$,

$$\xi_t(\lambda)\begin{pmatrix} 1 \\ \phi_t(\lambda,\nu) \end{pmatrix} = \begin{pmatrix} \nu \\ \nu\phi_t(\lambda,\nu) \end{pmatrix} ,$$

so that $\phi_t$ is indeed the projectivized eigenline of $\xi_t$. In particular, for all $\lambda$, the eigenvalues of $\xi(\lambda)$ are $\pm\nu$, so that

$$\text{Det}(\xi(\lambda)) = -\nu^2 = \frac{1}{\lambda}\Delta(\lambda) ,$$

and $\xi$ is thus a potential of $\Sigma$. This completes the proof. $\square$

It remains only to determine the derivative of this flow in the space of potentials. First, given a $2 \times 2$ matrix

$$M = \begin{pmatrix} a & b \\ c & d \end{pmatrix} ,$$

we denote

$$M_- := \begin{pmatrix} 0 & 0 \\ c & 0 \end{pmatrix} , \quad M_0 := \begin{pmatrix} a & 0 \\ 0 & d \end{pmatrix} , \quad \text{and } M_+ := \begin{pmatrix} 0 & b \\ 0 & 0 \end{pmatrix} ,$$

and, given any matrix-valued finite Laurent polynomial

$$N := \sum_k \lambda^k N_k ,$$

we denote

$$\Pi_-(N) := \sum_{k<0} \lambda^k N_k + N_{0,-} + \frac{1}{2}N_{0,0} , \text{ and}$$

$$\Pi_+(N) := \sum_{k>0} \lambda^k N_k + N_{0,+} + \frac{1}{2}N_{0,0} .$$
(3.6)

Given any potential $\xi$, we denote

$$\hat\xi[P,Q] := \Pi_-\left(P\left(\frac{1}{\lambda}\right)\xi(\lambda)\right) + \Pi_+\left(\frac{1}{\lambda^{g-1}}\overline{Q}(\lambda)\xi(\lambda)\right) .$$

**Lemma 3.2**

*The time derivative of $\xi_t$ is given by*

$$\partial_t \xi_t = [\hat\xi_t[P,Q], \xi_t] .$$
(3.7)

**Proof:** It suffices to work at $t=0$. Upon differentiating (3.1), we obtain, at $t=0$,

$$\partial_t f(\lambda,\nu) = \left(\nu P\left(\frac{1}{\lambda}\right) + \partial_t \text{Ln}(\rho)\right) f(\lambda,\nu) + O(\lambda\nu)\text{Exp}\left(t\nu P\left(\frac{1}{\lambda}\right)\right) ,$$

$$\partial_t g(\lambda,\nu) = \left(\nu P\left(\frac{1}{\lambda}\right) + \partial_t \text{Ln}(\sigma)\right) g(\lambda,\nu) + O(\lambda)\text{Exp}\left(t\nu P\left(\frac{1}{\lambda}\right)\right) ,$$

$$\partial_t (f \circ \sigma_3)(\lambda,\nu) = \left(\bar\nu \overline{Q}\left(\frac{1}{\bar\lambda}\right) - \partial_t \text{Ln}(\rho)\right)(f \circ \sigma_3)(\lambda,\nu) + O(\lambda\nu)\text{Exp}\left(t\bar\nu\overline{Q}\left(\frac{1}{\bar\lambda}\right)\right) , \text{ and}$$

$$\partial_t (g \circ \sigma_3)(\lambda,\nu) = \left(\bar\nu \overline{Q}\left(\frac{1}{\bar\lambda}\right) - \partial_t \text{Ln}(\sigma)\right)(g \circ \sigma_3)(\lambda,\nu) + O(1)\text{Exp}\left(t\bar\nu\overline{Q}\left(\frac{1}{\bar\lambda}\right)\right) .$$

Since $(f,g)^t$ is an eigenvector of $\xi$ with eigenvalue $\nu$, it follows that

$$\begin{pmatrix} \partial_t f(\lambda,\nu) \\ \partial_t g(\lambda,\nu) \end{pmatrix} - P\Big(\frac{1}{\lambda}\Big)\xi(\lambda)\begin{pmatrix} f(\lambda,\nu) \\ g(\lambda,\nu) \end{pmatrix} = \begin{pmatrix} \partial_t \mathrm{Ln}(\rho)f(\lambda,\nu) \\ \partial_t \mathrm{Ln}(\sigma)g(\lambda,\nu) \end{pmatrix} + \mathrm{Exp}\Big(t\nu P\Big(\frac{1}{\lambda}\Big)\Big)\begin{pmatrix} \mathrm{O}(\lambda\nu) \\ \mathrm{O}(\lambda) \end{pmatrix},$$

so that

$$(\partial_t - \hat{\xi}[P,Q])\begin{pmatrix} f(\lambda,\nu) \\ g(\lambda,\nu) \end{pmatrix} = \begin{pmatrix} (\partial_t\mathrm{Ln}(\rho) + (\omega_0^P/2) - (\omega_0^{\overline{Q}}/2))f(\lambda,\nu) \\ (\partial_t\mathrm{Ln}(\sigma) - (\omega_0^P/2) + (\omega_0^{\overline{Q}}/2))g(\lambda,\nu) \end{pmatrix} + \mathrm{Exp}\Big(t\nu P\Big(\frac{1}{\lambda}\Big)\Big)\begin{pmatrix} \mathrm{O}(\lambda\nu) \\ \mathrm{O}(\lambda) \end{pmatrix},$$

where here $\omega_0^P$ and $\omega_0^{\overline{Q}}$ denote respectively the upper left entries of the zero'th order terms of $P(1/\lambda)\xi(\lambda)$ and $\overline{Q}(1/\lambda)\lambda^{-(g-1)}\xi(\lambda)$. A similar reasoning yields

$$(\partial_t - \hat{\xi}[P,Q])\begin{pmatrix} (f\circ\sigma_3)(\lambda,\nu) \\ (g\circ\sigma_3)(\lambda,\nu) \end{pmatrix} = \begin{pmatrix} (-\partial_t\mathrm{Ln}(\rho) - (\omega_0^P/2) + (\omega_0^{\overline{Q}}/2))(f\circ\sigma_3)(\lambda,\nu) \\ (-\partial_t\mathrm{Ln}(\sigma) + (\omega_0^P/2) - (\omega_0^{\overline{Q}}/2))(g\circ\sigma_3)(\lambda,\nu) \end{pmatrix} + \mathrm{Exp}\Big(t\overline{\nu}\overline{Q}\Big(\frac{1}{\overline{\lambda}}\Big)\Big)\begin{pmatrix} \mathrm{O}(\lambda\nu) \\ \mathrm{O}(1) \end{pmatrix}.$$

Since $D$ is a non-special divisor, it follows that

$$(\partial_t - \hat{\xi}[P,Q])\begin{pmatrix} f \\ g \end{pmatrix} = \begin{pmatrix} af \\ bg \end{pmatrix},$$

for some constants $a$ and $b$. Comparing asymptotic expansions at $0$ and $\infty$ yields

$$\partial_t\mathrm{Ln}(\rho) = -\partial_t\mathrm{Ln}(\sigma) = -\frac{1}{2}\omega_0^P + \frac{1}{2}\omega_0^{\overline{Q}}, \tag{3.8}$$

so that $a = b = 0$, and

$$\begin{pmatrix} \partial_t f \\ \partial_t g \end{pmatrix} = \hat{\xi}[P,Q](\lambda)\begin{pmatrix} f \\ g \end{pmatrix}. \tag{3.9}$$

Since, for all $(\lambda,\nu)$, $(f_t(\lambda,\nu), g_t(\lambda,\nu))^t$ is an eigenvector of $\xi_t(\lambda)$ with eigenvalue $\nu$, it follows that

$$\partial_t \xi_t = [\hat{\xi}_t[P,Q], \xi_t],$$

as desired. $\square$

**3.2 - Complex solutions.** We now review how Baker–Akhiezer flows yield complex solutions of the sinh-Gordon equation. Let $z_i \in \mathrm{Jac}(\Sigma)$ be a potential point and, as in the introduction, define $\phi[z_0] : \mathbb{C} \to \mathrm{Jac}(\Sigma)$ by

$$\phi[z_0](w) := z_0 + \mathrm{V}[w,-w].$$

For $w$ sufficiently small, $\phi[z_0](w)$ is also a potential point, and we denote by $\xi(w)$ its potential, which we recall is well-defined up to multiplication by $-1$. This potential has the form

$$\xi(w) := \sum_{k=-1}^{g} \lambda^k \xi_k(w),$$

and, in particular,

$$\xi_{-1}(w) = \begin{pmatrix} 0 & \frac{1}{4i}e^{u(w)} \\ 0 & 0 \end{pmatrix}, \tag{3.10}$$

for some function $u$. In [3], Bobenko shows that this function is given explicitly by (1.3). We now show that it solves the sinh-Gordon equation (1.1). Indeed, by (3.4),

$$e^u = \frac{\rho}{\sigma},$$

where $\sigma$ and $\rho$ are as in (3.1). Using also (2.20), (2.21), (3.5), and (3.8), we then show that

$$\hat{\xi}_x := \hat{\xi}[1,-1] = \frac{1}{\lambda}\begin{pmatrix} 0 & \frac{1}{4i}e^u \\ 0 & 0 \end{pmatrix} + \begin{pmatrix} \frac{i}{2}u_y & \frac{1}{4i}e^{-u} \\ \frac{1}{4i}e^{-u} & -\frac{i}{2}u_y \end{pmatrix} + \lambda\begin{pmatrix} 0 & 0 \\ \frac{1}{4i}e^u & 0 \end{pmatrix},$$

and

$$\hat{\xi}_y := \hat{\xi}[i,-i] = \frac{1}{\lambda}\begin{pmatrix} 0 & \frac{1}{4}e^u \\ 0 & 0 \end{pmatrix} + \begin{pmatrix} -\frac{i}{2}u_x & -\frac{1}{4i}e^{-u} \\ \frac{1}{4i}e^{-u} & \frac{i}{2}u_x \end{pmatrix} + \lambda\begin{pmatrix} 0 & 0 \\ -\frac{1}{4i}e^u & 0 \end{pmatrix},$$

By (3.9), the pair $(\hat{\xi}_x, \hat{\xi}_y)$ satisfies the Lax equation

$$\partial_y \hat{\xi}_x - \partial_x \hat{\xi}_y + [\hat{\xi}_x, \hat{\xi}_y] = 0,$$

from which it then follows that

$$u_{z\bar{z}} + \frac{1}{8}\sinh(2u) = 0,$$

as desired.

**3.3 - The real locus.** A potential of $\Sigma$ corresponds to a real solution of the sinh-Gordon equation whenever it satisfies the *reality condition*

$$\xi(\lambda) = \lambda^{g-1}\overline{\xi\left(\frac{1}{\bar{\lambda}}\right)}^t. \qquad (3.11)$$

Indeed (2.2), (2.20) and (3.11) together imply that

$$\sigma_{-1}^2 = -|\sigma_{-1}|^2,$$

so that

$$\sigma_{-1} \in i\mathbb{R} \setminus \{0\}, \qquad (3.12)$$

which, by (3.10), indeed guarantees reality of our solutions. Note that the space of real potentials divides into two parts, determined by the sign of $i\sigma_{-1}$, and we henceforth consider only those real potentials satisfying

$$i\sigma_{-1} > 0.$$

By Remark 2.3, this is always possible and, furthermore, resolves the order 2 ambiguity highlighted there.

We now show that the set of points of the Jacobi variety which are divisors of potentials satisfying (3.11) is a real $g$-dimensional subspace of $\mathrm{Jac}(\Sigma)$, given by the vanishing of a certain affine condition. Indeed, we define $\hat{\Phi} : \mathbb{C}^g \to \mathbb{C}^g$ by

$$\hat{\Phi}(z)_i := z_i - \overline{z}_{(g+1)-i}.$$

Using (2.9) and (2.10), it is straightforward to show that $\hat{\Phi}$ sends periods to periods, and thus descends to a function $\Phi : \mathrm{Jac}(\Sigma) \to \mathrm{Jac}(\Sigma)$.

**Theorem & Definition 3.3**

*A point $z \in \mathrm{Jac}(\Sigma)$ is the divisor of a potential satisfying the reality condition (3.11) if and only if*

$$\Phi(z) = 0. \qquad (3.13)$$

We call the set of all such points the *real locus of* $\mathrm{Jac}(\Sigma)$, *and we denote it by* $\mathrm{Jac}_{\mathrm{re}}(\Sigma)$.

**Remark 3.1.** In particular, every point of the real locus is the image under the Abel map of a non-special divisor supported on $\Sigma^*$.

Theorem 3.3 follows immediately from the following two lemmas. We first prove the result for potential points, and we then show that every point of $\mathrm{Jac}_{\mathrm{re}}(\Sigma)$ is a potential point.

**Lemma 3.4**

If $z \in \mathrm{Jac}(\Sigma)$ is the image under the Abel map of a non-special degree $g$ divisor supported on $\Sigma^*$, then its potential $\xi$ satisfies the reality condition (3.11) if and only if,

$$\Phi(z) = 0 . \tag{3.14}$$

**Proof:** Indeed, denote

$$\xi := \begin{pmatrix} \omega & \sigma \\ \tau & -\omega \end{pmatrix} ,$$

and suppose that this potential satisfies the reality condition (3.11). Observe that the projectivized eigenlines of $\xi$ and $\lambda^{g-1}\overline{\xi(1/\bar\lambda)}^t$ are given respectively by

$$\phi(\lambda, \nu) = \frac{\nu - \omega(\lambda)}{\sigma(\lambda)} , \text{ and}$$

$$\psi(\lambda, \nu) = \frac{\nu - \lambda^{g-1}\overline{\omega}(1/\lambda)}{\lambda^{g-1}\overline{\tau}(1/\lambda)} .$$

By hypothesis, these two functions coincide. We verify that their respective pole divisors satisfy

$$[\phi]_\infty = [\sigma \circ \pi]_0 - [\nu - \omega \circ \pi]_0 + (\infty) , \text{ and}$$
$$[\overline{\psi} \circ \sigma_3 \circ \sigma_1]_\infty = [\tau \circ \pi]_0 - [\nu + \omega \circ \pi]_0 + (0) .$$

However, by (2.22), since $(\nu \pm \omega \circ \pi)$ and $\sigma$ all have poles at $0$ and $\infty$ and only at these points, and since $\tau$ is non-vanishing at $0$, and has a pole only at $\infty$,

$$[\nu - \omega \circ \pi]_0 + [\nu + \omega \circ \pi]_0 = [\sigma \circ \pi]_0 + [\tau \circ \pi]_0 ,$$

so that

$$[\overline{\psi} \circ \sigma_3 \circ \sigma_1]_\infty = [\nu - \omega \circ \pi]_0 - [\sigma \circ \pi]_0 + (0) = -[\phi]_\infty + (\infty) + (0) .$$

Since $\phi = \psi$, the finite part of its pole divisor thus satisfies

$$[\phi]_{\infty,\mathrm{fin}} = [\phi]_\infty - (\infty) = [\psi]_\infty - (\infty) = -(\sigma_1 \circ \sigma_3)([\phi]_\infty) + (0) = -(\sigma_1 \circ \sigma_3)([\phi]_{\infty,\mathrm{fin}}) .$$

Applying the Abel map to both sides of this identity yields

$$\mathcal{A}([\phi]_{\infty,\mathrm{fin}}) + \mathcal{A}((\sigma_1 \circ \sigma_3)([\phi]_{\infty,\mathrm{fin}})) = 0 .$$

By (2.12), for all $p \in \Sigma$,

$$\mathcal{A}((\sigma_1\sigma_3)(p))_i = -\mathcal{A}(\infty)_i - \overline{\mathcal{A}(p)}_{(g+1)-i} ,$$

so that

$$\mathcal{A}([\phi]_{\infty,\mathrm{fin}})_i - \overline{\mathcal{A}([\phi]_{\infty,\mathrm{fin}})}_{(g+1)-i} = g\mathcal{A}(\infty)_i .$$

Finally, since $(0)^2 - (\infty)^2$ is the divisor of the meromorphic function $\lambda$, and since $\mathcal{A}(0) = 0$,

$$2\mathcal{A}(\infty) = 2\mathcal{A}(\infty) - 2\mathcal{A}(0) = 0 , \tag{3.15}$$

so that, since the genus $g$ is even,

$$\Phi^\varepsilon(z)_i = \mathcal{A}([\phi]_{\infty,\mathrm{fin}})_i - \overline{\mathcal{A}([\phi]_{\infty,\mathrm{fin}})}_{(g+1)-i} = 0 ,$$

as desired. The converse follows upon applying the same argument in the reverse direction, and this completes the proof. $\square$

**Lemma 3.5**

*Every point of the real locus lies in the image under the Abel map of the set of non-special degree $g$ divisors supported on $\Sigma^*$.*

**Remark 3.2.** In particular, all real, finite-type solutions of (1.1) are regular over the whole of $\mathbb{R}^2$. Indeed, using the notation of (1.3), the point $(x,y)$ is a singular point of the solution $u[z_0]$ if and only if the point $\phi[z_0](x+iy, -(x+iy))$ is not a potential point of $\text{Jac}(\Sigma)$. Since every point of $\text{Jac}_{\text{re}}(\Sigma)$ is a potential point, regularity follows.

**Proof:** Let $D_g$ denote the set of non-special degree $g$ divisors over $\Sigma$, and let $D_g^*$ denote the subset consisting of those divisors that are supported over $\Sigma^*$. $\mathcal{A}(D_g^*) \cap \text{Jac}_{\text{re}}(\Sigma)$ is trivially an open subset of $\text{Jac}_{\text{re}}(\Sigma)$. We now claim that it closed. To show this, let $(z_m)_{m \in \mathbb{N}}$ be a sequence of points of this set converging to some limit $z_\infty$, say. For all $m$, let $D_m$ denote the preimage of $z_m$ under the Abel map, and let $\xi_m$ denote its associated potential. By Lemma 2.14, it will suffice to show that $(\xi_m)_{m \in \mathbb{N}}$ is precompact. However, for all unit $\lambda$, the matrix $\lambda^{-(g-1)/2} i^{(1-p)} \xi(\lambda)$ is trace-free and skew-adjoint, so that

$$\|\xi(\lambda)^2\| = \text{Det}(\lambda^{-(g-1)/2} i^{(1-p)} \xi(\lambda)) = \lambda^{-(g-1)} \varepsilon^{(1-p)} \text{Det}(\xi(\lambda)) = \varepsilon^{(1-p)} \lambda^{-g} \Delta(\lambda) \ .$$

Thus, bearing in mind the Cauchy-Schwarz inequality, for each $i$,

$$\|\xi_i\|^2 = \left\| \frac{1}{2\pi i} \int_{|\lambda|=1} \frac{1}{\lambda^{i+1}} \xi(\lambda) d\lambda \right\|^2 \leqslant \frac{1}{2\pi} \int_{|\lambda|=1} \|\xi(\lambda)\|^2 d\lambda = \frac{1}{2\pi} \int_{|\lambda|=1} \varepsilon^{(1-p)} \lambda^{-g} \Delta(\lambda) d\lambda \leqslant \sup_{|\lambda|=1} |\Delta(\lambda)| \ .$$

The sequence $(\xi_m)_{m \in \mathbb{N}}$ is thus uniformly bounded. By the Heine-Borel theorem, it is therefore precompact, and it follows that $\mathcal{A}(D_g^*) \cap \text{Jac}_{\text{re}}(\Sigma)$ is a closed subset of $\text{Jac}_{\text{re}}(\Sigma)$, as asserted.

We now show that $\mathcal{A}(D_g^*) \cap \text{Jac}_{\text{re}}(\Sigma)$ is a dense subset of $\text{Jac}_{\text{re}}(\Sigma)$. Indeed, suppose the contrary, and choose $z \in \text{Jac}_{\text{re}}(\Sigma)$ such that there exists a neighbourhood $\Omega$, say, of this point in $\text{Jac}_{\text{re}}(\Sigma)$ not lying in $\mathcal{A}(D_g)$. Recall that the complement of $\mathcal{A}(D_g)$ is everywhere locally given by the zero set of some holomorphic function. However, since $\text{Jac}_{\text{re}}(\Sigma)$ is a real submanifold of maximal dimension, any holomorphic function vanishing over $\Omega$ vanishes identically. This is absurd, since $\mathcal{A}(D_g)$ is non-empty, and we conclude that $\mathcal{A}(D_g) \cap \text{Jac}_{\text{re}}(\Sigma)$ is dense in $\text{Jac}_{\text{re}}(\Sigma)$. A similar argument also shows that $\mathcal{A}(D_g^*) \cap \text{Jac}_{\text{re}}(\Sigma)$ is dense in $\text{Jac}_{\text{re}}(\Sigma)$, and this completes the proof. $\square$

## 4 - Boundary conditions.

**4.1 - The Durham locus.** We now turn our attention to expressing the Sklyanin condition in terms of the Jacobi variety. First, for a non-negative integer $q$, we say that a potential $\xi$ satisfies the *$q$-Sklyanin condition* whenever

$$K(\lambda)\xi(\lambda) = \lambda^q \xi\left(\frac{1}{\lambda}\right) K(\lambda) \ , \tag{4.1}$$

where $K(\lambda)$ denotes Sklyanin's $K$-matrix

$$K(\lambda) := \begin{pmatrix} 4A - 4B\lambda & \lambda - 1/\lambda \\ \lambda - 1/\lambda & 4A - 4B/\lambda \end{pmatrix} \ ,$$

and where $A$ and $B$ are real constants. For ease of presentation, we will only address in detail the case where $|A| \neq |B|$. The modifications required to address the simpler cases where $|A| = |B| \neq 0$ and $A = B = 0$ will be explained at the end of this section.

In [19], we show that the solution $u$ given by (3.10) satisfies the Durham condition (1.2) along the line $\mathbb{R} \times \{0\}$ if and only if $\xi$ satisfies the $(g-1)$-Sklyanin condition. Note, however, that (4.1) differs from Condition (3.3) of [19]. This is because we are only concerned in [19] with those potentials which also satisfy the reality condition (3.11), for which (4.1) and Condition (3,3) of [19] are equivalent. However, without the reality condition, (4.1) strikes us as more natural since, as we will see in Section 4.2, it is preserved under a large family of flows, whilst it is not clear that the same holds true for Condition (3.3) of [19].

We define the *Durham locus* $\mathrm{Dur}(\Sigma) \subseteq \mathrm{Jac}(\Sigma)$ to be the set of all images under the Jacobi map of divisors of potentials satisfying the $(g-1)$-Sklyanin condition. The purpose of this and the following section is to describe the geometry of the Durham locus, from which Theorems 1.1 and 1.2 will follows. We will see that it consists of at most 2 connected components, each of which is a Zariski open subset of a $(g/2)$-dimensional complex affine subspace of $\mathrm{Jac}(\Sigma)$.

We first identify exceptional points of the spectral curve $\Sigma$. We define the *Sklyanin subset* $S \subseteq \Sigma$ by

$$S := \{(\lambda, \nu) \in \Sigma \mid \mathrm{Det}(K(\lambda)) = 0\} ,$$

and we call its elements *Sklyanin points*.

**Lemma 4.1**

*The Sklyanin subset has cardinality $8$ and is invariant under the actions of the involutions $\sigma_1$, $\sigma_2$ and $\sigma_3$.*

**Proof:** Indeed, the determinant of the $K$-matrix is given by

$$\mathrm{Det}(K(\lambda)) = 16A^2 + 16B^2 - 16AB\left(\lambda + \frac{1}{\lambda}\right) - \left(\lambda - \frac{1}{\lambda}\right)^2 .$$

In particular, $\lambda^2 \mathrm{Det}(K(\lambda))$ is a polynomial of order $4$ with real coefficients, so that $\mathrm{Det}(K(\lambda))$ has $4$ roots over $\hat{\mathbb{C}}$. Since

$$\lim_{\lambda \to 0} \mathrm{Det}(K(\lambda)) = \lim_{\lambda \to \infty} \mathrm{Det}(K(\lambda)) = -\infty ,$$

and

$$\mathrm{Det}(K(\pm 1)) = 16(A \mp B)^2 > 0 ,$$

it follows by the intermediate value theorem that $\mathrm{Det}(K(\lambda))$ has $4$ distinct roots over the real line. Since $\Sigma$ has no branch points over the real line, it follows that the Sklyanin subset has cardinality $8$, as desired. The Sklyanin subset is trivially invariant under the action of $\sigma_1$. Since $\mathrm{Det}(K(\lambda))$ is real, it is invariant under the action of $\sigma_2$. Finally, since $\mathrm{Det}(K(\lambda)) = \mathrm{Det}(K(1/\lambda))$, it is invariant under the action of $(\sigma_2 \circ \sigma_3)$, and therefore also $\sigma_3$. This completes the proof. $\square$

We now define $\hat{\Psi} : \mathbb{C}^g \to \mathbb{C}^g$ by

$$\hat{\Psi}(z)_i := z_i - z_{(g+1)-i} .$$

As before, it is straightforward to show that $\hat{\Psi}$ maps periods to periods and thus descends to a function $\Psi : \mathrm{Jac}(\Sigma) \to \mathrm{Jac}(\Sigma)$. We say that a subset $S_0$ of $S$ is a *special Sklyanin subset* whenever it has cardinality $4$ and

$$S = S_0 \sqcup (\sigma_2 \sigma_3)(S_0) .$$

Note that there are precisely $4$ such subsets. Given such a special Sklyanin subset $S_0$, we define its *Sklyanin subspace* $\mathrm{Skl}(S_0)$ by

$$\mathrm{Skl}(S_0) := \{z \in \mathrm{Jac}(\Sigma) \mid \Psi(z) = \mathcal{A}(S_0)\} .$$

There are at most $4$ Sklyanin subspaces, each of which is a complex $(g/2)$-dimensional affine subspace of the Jacobi variety.

We will show in this section that the Durham locus of $\Sigma$ is contained in the union of the Sklyanin subspaces, and we will show in the next section that every connected component is a Zariski open subset of some Sklyanin subspace. Let $\xi$ be a potential of $\Sigma$ satisfying the $(g-1)$-Sklyanin condition, and let

$$\phi(\lambda, \nu) := \frac{\nu - \omega(\lambda)}{\sigma(\lambda)}$$

denote its projectivized eigenline.

**Lemma 4.2**

*At every root $\lambda$ of $\mathrm{Det}(K(\lambda))$, the kernel of $K(\lambda)$ is an eigendirection of $\xi(\lambda)$.*

**Proof:** Indeed, suppose the contrary, and let $u$ be a non-trivial element of $\mathrm{Ker}(K(\lambda))$. Since $u$ is not an eigenvector of $\xi(\lambda)$, $\xi(\lambda)u \notin \mathrm{Ker}(K(\lambda))$, and so

$$0 = \lambda^{g-1}\xi\left(\frac{1}{\lambda}\right)K(\lambda)u = K(\lambda)\xi(\lambda)u \neq 0 \ .$$

This is absurd, and the result follows. $\square$

We say that a Sklyanin point $(\lambda,\nu)$ is a *special Sklyanin point* of $\xi$ whenever $\nu$ is the eigenvalue of $\mathrm{Ker}(K(\lambda))$. We denote the set of special Sklyanin points of $\xi$ by $S_0(\xi)$.

**Lemma 4.3**

*If $(\lambda,\nu)$ is a special Sklyanin point of $\xi$, then $(\sigma_1\sigma_2\sigma_3)(\lambda,\nu) = (1/\lambda, -\nu/\lambda^{g-1})$ is also a special Sklyanin point of $\xi$. In particular*

$$S = S_0(\xi) \cup (\sigma_2\sigma_3)(S_0(\xi)) \ , \tag{4.2}$$

*so that $S_0(\xi)$ is a special Sklyanin subset.*

**Proof:** Let $u$ be an eigenvector of $\xi(\lambda)$ with eigenvalue $-\nu$. In particular, $K(\lambda)u \neq 0$, and

$$\lambda^{g-1}\xi\left(\frac{1}{\lambda}\right)K(\lambda)u = \varepsilon K(\lambda)\xi(\lambda)u = -\nu K(\lambda)u \ .$$

It follows that $K(\lambda)u$ is an eigenvector of $\xi(1/\lambda)$ with eigenvalue $-\nu/\lambda^{-(g-1)}$. However,

$$K\left(\frac{1}{\lambda}\right)K(\lambda)u = \mathrm{Det}(K(\lambda))u = 0 \ .$$

It follows that $K(\lambda)u$ is a non-trivial element of $\mathrm{Ker}(K(1/\lambda))$, and the result follows. $\square$

**Lemma 4.4**

*If $\xi$ satisfies the $(g-1)$-Sklyanin condition, then, for every special Sklyanin point $(\lambda,\nu)$,*

$$\phi(\lambda,\nu) = -\frac{(4A - 4B\lambda)}{(\lambda - 1/\lambda)} \ . \tag{4.3}$$

**Remark 4.1.** Note that when $B$ vanishes, this value is equal to $\pm 1$.

**Proof:** Indeed, note first that, since $|A| \neq |B|$, $(\pm 1)$ is not a root of $\mathrm{Det}(\mathrm{Ker}(\lambda))$, and the right-hand side of (4.3) is finite. We verify by inspection that when $\mathrm{Det}(K(\lambda)) = 0$,

$$\mathrm{Ker}(K(\lambda)) = \langle (1, -(4A - 4B\lambda)/(\lambda - 1/\lambda))^t \rangle \ ,$$

and the result follows by Lemma 4.2. $\square$

**Theorem 4.5**

*If $\xi$ satisfies the $(g-1)$-Sklyanin condition, then*

$$\Psi(\mathcal{A}([\phi]_{\infty,\mathrm{fin}})) = \mathcal{A}(S_0) \ . \tag{4.4}$$

*Conversely, if $\phi$ satisfies (4.4) and if, in addition, it satisfies (4.3) at every special Sklyanin point, then $\xi$ satisfies the $(g-1)$-Sklyanin condition.*

In order to prove Theorem 4.5, consider the function

$$\psi(\lambda,\nu) := \frac{(\lambda - 1/\lambda) + (4A - 4B/\lambda)\phi(\lambda,\nu)}{(\lambda - 1/\lambda)\phi(\lambda,\nu) + (4A - 4B\lambda)} \ . \tag{4.5}$$

Note that $\xi$ satisfies the Sklyanin condition if and only if $\psi$ is the projectivized eigenline of $\xi(1/\lambda)$, that is, if and only if

$$\psi(\lambda,\nu) = \phi\left(\frac{1}{\lambda}, \frac{\nu}{\lambda^{g-1}}\right) \ . \tag{4.6}$$

Let $N_\psi$ and $D_\psi$ denote respectively the numerator and denominator of (4.5).

**Lemma 4.6**

$N_\phi$ and $D_\phi$ can only both vanish at the special Sklyanin points of $\xi$.

**Proof:** Indeed, $N_\psi$ and $D_\psi$ both vanish at $(\lambda, \nu)$ if and only if $(1, \phi(\lambda, \nu))^t$ is an element of $\mathrm{Ker}(K(\lambda))$, that is, if and only if $(\lambda, \nu)$ is a special Sklyanin point of $\xi$, as desired. $\square$

**Lemma 4.7**

If $\phi$ satisfies (4.3), then, at every special Sklyanin point $(\lambda, \nu)$ of $\xi$,

$$\mathrm{Min}(\mathrm{Ord}(N_\psi(\lambda, \nu)), \mathrm{Ord}(D_\psi(\lambda, \nu))) = 1 \ . \tag{4.7}$$

**Proof:** Indeed, let $m \geqslant 1$ denote the minimum of the orders of these two zeroes. Since

$$(4A - 4B\lambda)(4A - 4B/\lambda) = (\lambda - 1/\lambda)^2 \neq 0 \ ,$$

neither of the factors on the left hand side vanishes at $(\lambda, \nu)$, and $\phi$ therefore has neither a zero nor a pole at this point. Let $\zeta$ be a local coordinate of $\Sigma$ about $(\lambda, \nu)$. We have

$$(\lambda - 1/\lambda) + (4A - 4B/\lambda)\phi(\lambda, \nu) = \mathrm{O}(\zeta^m) \ , \text{ and}$$
$$(4A - 4B\lambda) + (\lambda - 1/\lambda)\phi(\lambda, \nu) = \mathrm{O}(\zeta^m) \ .$$

Thus

$$\frac{(4A - 4B\lambda)}{(\lambda - 1/\lambda)} + \mathrm{O}(\zeta^m) = \phi(\lambda, \nu) = \frac{(\lambda - 1/\lambda)}{(4A - 4B/\lambda)} + \mathrm{O}(\zeta^m) \ ,$$

so that

$$\mathrm{Det}(K(\lambda)) = (4A - 4B\lambda)(4A - 4B/\lambda) - (\lambda - 1/\lambda)^2 = \mathrm{O}(\zeta^m) \ .$$

Since $\mathrm{Det}(K(\lambda))$ only has simple zeroes, and since $\Sigma$ is a graph over the first component at $(\lambda, \nu)$, it follows that $m = 1$, as desired. $\square$

**Lemma 4.8**

If $\phi$ satisfies (4.3) at each Sklyanin point, then

$$\mathcal{A}([\psi]_{\infty,\mathrm{fin}}) = \mathcal{A}([\phi]_{\infty,\mathrm{fin}}) + \mathcal{A}(\infty) - \mathcal{A}(S_0) \ . \tag{4.8}$$

**Proof:** Denote $P := \pi^{-1}(\{\pm 1\})$ and note that, since $|A| \neq |B|$, $P$ and $S$ are disjoint. Consider first the case where no point of $P$ is a pole of $\phi$. In particular, since every pole of $\phi$ then appears with non-trivial coefficient in $D_\psi$, none of these poles contributes a pole to $\psi$. It follows that the poles of $\psi$ over $\Sigma^*$ are precisely those zeroes of the denominator over this set which are not cancelled by zeroes of the numerator. By Lemma 4.6, zeroes of the denominator are cancelled by zeroes of the numerator precisely at the special Sklyanin points of $\xi$. Furthermore, by Lemma 4.7, precisely one zero of the denominator is cancelled by a zero of the numerator at every such point. Denoting by $Q$ the greatest divisor bounded above by both $[\phi]_0$ and $[4A - 4B\lambda]_0$, we therefore obtain

$$\begin{aligned}
\mathcal{A}([\psi]_{\infty,\mathrm{fin}}) &= \mathcal{A}([(\lambda - 1/\lambda)\phi + (4A - 4B\lambda)]_{0,\mathrm{fin}}) - \mathcal{A}(S_0) \\
&= \mathcal{A}([(\lambda - 1/\lambda)\phi/(4A - 4B\lambda) + 1]_{0,\mathrm{fin}}) + \mathcal{A}(Q) - \mathcal{A}(S_0) \\
&= \mathcal{A}([(\lambda - 1/\lambda)\phi/(4A - 4B\lambda)]_0) + \mathcal{A}(Q) - \mathcal{A}(S_0) \\
&= \mathcal{A}([\phi]_0) + \mathcal{A}(P) - \mathcal{A}(S_0) \ .
\end{aligned}$$

Since $(0)^4 - P$ is the divisor of the meromorphic function $\lambda^2/(\lambda^2 - 1)$, and since $\mathcal{A}(0) = 0$,

$$\mathcal{A}(P) = 4\mathcal{A}(0) = 0 \ .$$

Furthermore,
$$\mathcal{A}([\phi]_0) = \mathcal{A}([\phi]_\infty) = \mathcal{A}([\phi]_{\infty,\text{fin}}) + \mathcal{A}(\infty) ,$$

and the result follows in the case where $\phi$ has no pole over $\pm 1$ upon combining these identities.

Suppose now that $\phi$ has at least one pole over $\pm 1$. In order to understand the general case, it is sufficient to address the case where $\phi$ has a simple pole at the point $z_0 := (\lambda_0, \nu_0)$, say. In this case, $\psi$ has a simple pole at this point. On the other hand, since the pole of $\phi$ at this point is eliminated by the factor $(\lambda - 1/\lambda)$, the function $(\lambda - 1/\lambda)\phi + (4A - 4B\lambda)$ has a total of $(g+3)$ poles over $\Sigma$. It therefore has as many zeroes over this surface, and since 4 of these zeroes are accounted for by the special Sklyanin points, it has precisely $(g-1)$ remaining zeroes, defining a divisor $D'$ which yields the remaining finite poles of $\psi$. Consequently, as before,
$$\mathcal{A}([\psi]_{\infty,\text{fin}}) = \mathcal{A}((z_0)) + \mathcal{A}([(\lambda - 1/\lambda)\phi + (4A - 4B\lambda)]_{0,\text{fin}}) - \mathcal{A}(S_0) .$$

However, since the zero of $(\lambda - 1/\lambda)$ at $z_0$ is eliminated by the pole of $\phi$ at this point,
$$\mathcal{A}([(\lambda - 1/\lambda)\phi/(4A - 4B\lambda)]_0) = \mathcal{A}([\phi]_0) + \mathcal{A}(P) - \mathcal{A}((z_0)) - \mathcal{A}(Q) ,$$

and the result now follows as before. $\square$

**Proof of Theorem 4.5:** Indeed, suppose first that $\xi$ satisfies the $(g-1)$-Sklyanin condition and let $S_0$ denote its special Sklyanin subset. By (4.6) and (4.8),

$$\begin{aligned}\mathcal{A}([\phi]_{\infty,\text{fin}}) &= \mathcal{A}((\sigma_2\sigma_3)([\psi]_{\infty,\text{fin}})) \\ &= \mathcal{A}((\sigma_2\sigma_3)([\phi]_{\infty,\text{fin}})) + \mathcal{A}((\sigma_2\sigma_3)(\infty)) - \mathcal{A}((\sigma_2\sigma_3)(S_0)) \\ &= \mathcal{A}((\sigma_2\sigma_3)([\phi]_{\infty,\text{fin}})) + \mathcal{A}(0) - \mathcal{A}(S) + \mathcal{A}(S_0) .\end{aligned}$$

Since $\mathcal{A}(0) = 0$, and since $S - (0)^8$ is the divisor of the meromorphic function $\lambda^{-2}\text{Det}(K(\lambda))$,
$$\mathcal{A}([\phi]_{\infty,\text{fin}}) = \mathcal{A}((\sigma_2\sigma_3)([\phi]_{\infty,\text{fin}})) + \mathcal{A}(S_0) .$$

Observing that $\mathcal{A}(\infty)$ is real, by (2.12), for all $x \in \Sigma$, and for all $i$,
$$\mathcal{A}((\sigma_2\sigma_3)(x))_i = \mathcal{A}(\infty)_i + \mathcal{A}(x)_{(g+1)-i} .$$

Thus, for all $i$,
$$\mathcal{A}([\phi]_{\infty,\text{fin}})_i = g\mathcal{A}(\infty)_i + \mathcal{A}([\phi]_{\infty,\text{fin}})_{(g+1)-i} + \mathcal{A}(S_0)_i .$$

Since the genus $g$ is even, and since $(0)^2 - (\infty)^2$ is the divisor of the meromorphic function $\lambda$, for all $i$,
$$\begin{aligned}\mathcal{A}([\phi]_{\infty,\text{fin}})_i &= g\mathcal{A}(0)_i + \mathcal{A}([\phi]_{\infty,\text{fin}})_{(g+1)-i} + \mathcal{A}(S_0)_i \\ &= \mathcal{A}([\phi]_{\infty,\text{fin}})_{(g+1)-i} + \mathcal{A}(S_0)_i ,\end{aligned}$$

so that
$$\Psi(\mathcal{A}([\phi]_{\infty,\text{fin}})) = \mathcal{A}(S_0) ,$$

as desired.

Conversely, suppose that $\phi$ satisfies (4.3). If $\phi$ satisfies (4.4), then, by (4.8), the image under the Abel map of the finite part of its pole divisor coincides with the image under the Abel map of the finite part of the pole divisor of $\psi$. Since these divisors are non-special, it follows that $\phi = \psi$, from which it readily follows that $\xi$ satisfies the $(g-1)$-Sklyanin condition. This completes the proof. $\square$

We conclude this section by discussing the modifications required to address the case where $|A| = |B|$. Consider first the case where $A = B \neq 0$. In this case, the $K$-matrix factorizes as
$$K(\lambda) = \left(1 - \frac{1}{\lambda}\right) K'(\lambda) ,$$

where

$$K'(\lambda) := \begin{pmatrix} -4A\lambda & \lambda+1 \\ \lambda+1 & 4A \end{pmatrix} .$$

It now suffices to repeat the theory developed above with $K'$ instead of $K$. Note that there are now only 4 Sklyanin points, 2 special Sklyanin subsets, and 2 Sklyanin subspaces. The case where $A = -B \neq 0$ is likewise addressed in a similar manner. Finally, the case where $A = B = 0$ is addressed by replacing the $K$-matrix with

$$K''(\lambda) := \begin{pmatrix} 0 & 1 \\ 1 & 0 \end{pmatrix} .$$

This corresponds to the Neumann boundary conditions addressed by Bobenko-Kuksin in [4].

**4.2 - Preservation of the Sklyanin condition.** A priori, the criteria of Lemma 4.4 define a submanifold of the Jacobi variety of complex codimension 4. We now show that each Sklyanin subspace is either wholly contained in that submanifold, or only meets it at special points.

**Theorem 4.9**

*If the criteria of Lemma 4.4 are satisfied at a single non-special point of some Sklyanin subspace Skl($S_0$), then they are satisfied at every non-special point of this Sklyanin subspace.*

*In particular, since all real points are non-special, if these criteria are satisfied at a single real point of some Sklyanin subspace, then they are satisfied at every real point of that Sklyanin subspace.*

For all non-negative integer $k$, we define the projection $\Theta_k^q$ by

$$\Theta_k^q(\xi) = \Pi_-(\lambda^{-k}\xi) + \Pi_+(\lambda^{k-q}\xi) ,$$

where $\Pi_-$ and $\Pi_+$ are as in (3.6). Theorem 4.9 will be a consequence of the following result concerning the action of $\Theta_k^q$ on potentials satisfying the $q$-Sklyanin condition.

**Lemma 4.10**

*If $\xi$ satisfies the $q$-Sklyanin condition, then, for all non-negative integer $k$, $\Theta_k^q(\xi)$ satisfies the $(0,1)$-Sklyanin condition.*

Lemma 4.10 will follow from Lemmas 4.12 and 4.13. We first introduce a new variable $\gamma$ and apply the guage transform

$$\xi(\lambda, \gamma) := e^{-\theta\sigma_0/2} \xi\left(\frac{\lambda}{\gamma}\right) e^{\theta\sigma_0/2} , \tag{4.9}$$

where

$$e^\theta = \gamma . \tag{4.10}$$

This transformation clarifies symmetries of the $K$-matrix. Indeed,

$$e^{\theta\sigma_0/2} K\left(\frac{\lambda}{\gamma}\right) e^{\theta\sigma_0/2} = K(\lambda, \gamma) , \tag{4.11}$$

where

$$K(\lambda, \gamma) := \left(\frac{\lambda}{\gamma} - \frac{\gamma}{\lambda}\right) \begin{pmatrix} 0 & 1 \\ 1 & 0 \end{pmatrix} + \begin{pmatrix} 4A\gamma - 4B\lambda & 0 \\ 0 & 4A/\gamma - 4B/\lambda \end{pmatrix} . \tag{4.12}$$

Furthermore,

$$K(\lambda)\xi(\lambda) = \lambda^q \xi\left(\frac{1}{\lambda}\right) K(\lambda)$$

if and only if

$$K(\lambda, \gamma)\xi(\lambda, \gamma) = \lambda^q \gamma^{-q} \xi\left(\frac{1}{\lambda}, \frac{1}{\gamma}\right) K(\lambda, \gamma) . \tag{4.13}$$

Note that $\xi$ is a finite sum of the form

$$\xi(\lambda, \gamma) := \sum_{m,n} \xi_{m,n} \lambda^m \gamma^n. \tag{4.14}$$

where $\xi_{m,n}$ is non-zero only if $|m+n| \leqslant 1$. Furthermore, terms along the diagonal $m+n=1$ are multiples of $\sigma_-$ and only involve $\tau_*$, terms along the diagonal $m+n=0$ are multiples of $\sigma_0$ and only involve $\omega_*$, and terms along the diagonal $m+n=-1$ are multiples of $\sigma_+$ and only involve $\sigma_*$. Consequently, in this guage

$$\begin{aligned}
\Theta_k^q(\xi) = & \sum_{\substack{m+n=0 \\ m<0}} \lambda^m \gamma^n \omega_{m+k,n-k} \sigma_0 + \sum_{\substack{m+n=1 \\ m \leqslant 0}} \lambda^m \gamma^n \tau_{m+k,n-k} \sigma_- \\
& + \sum_{\substack{m+n=-1 \\ m<0}} \lambda^m \gamma^n \sigma_{m+k,n-k} \sigma_+ + \sum_{\substack{m+n=0 \\ m>0}} \lambda^m \gamma^n \omega_{m+q-k,n-q+k} \sigma_0 \\
& + \sum_{\substack{m+n=1 \\ m>0}} \lambda^m \gamma^n \tau_{m+q-k,n-q+k} \sigma_- + \sum_{\substack{m+n=-1 \\ m \geqslant 0}} \lambda^m \gamma^n \sigma_{m+q-k,n-q+k} \sigma_+ \\
& + \frac{1}{2} \big(\omega_{k,-k} + \omega_{q-k,-q+k}\big) \sigma_0 \ .
\end{aligned} \tag{4.15}$$

Substituting (4.14) into (4.13), we determine that the $q$-Sklyanin condition is equivalent to the condition that, for all $m,n$,

$$A_1^q(\xi; m, n) = A_2^q(\xi; m, n) = A_3^q(\xi; m, n) = 0, \tag{4.16}$$

where

$$\begin{aligned}
A_1^q(\xi; m, n) := & \tau_{m-1,n+1} - 4B\omega_{m-1,n} + 4A\omega_{m,n-1} - \tau_{m+1,n-1} \\
& + \sigma_{q-m-1,-q-n+1} - 4A\omega_{q-m,-q-n+1} \\
& + 4B\omega_{q-m+1,-q-n} - \sigma_{q-m+1,-q-n-1},
\end{aligned} \tag{4.17}$$

$$\begin{aligned}
A_2^q(\xi; m, n) := & -\omega_{m-1,n+1} - 4B\sigma_{m-1,n} + 4A\sigma_{m,n-1} + \omega_{m+1,n-1} \\
& + \omega_{q-m-1,-q-n+1} + 4B\sigma_{q-m-1,-q-n} \\
& - 4A\sigma_{q-m,-q-n-1} - \omega_{q-m+1,-q-n-1}, \text{ and}
\end{aligned} \tag{4.18}$$

$$\begin{aligned}
A_3^q(\xi; m, n) := & \omega_{m-1,n+1} + 4A\tau_{m,n+1} - 4B\tau_{m+1,n} - \omega_{m+1,n-1} \\
& - \omega_{q-m-1,-q-n+1} - 4A\tau_{q-m,-q-n+1} \\
& + 4B\tau_{q-m+1,-q-n} + \omega_{q-m+1,-q-n-1}.
\end{aligned} \tag{4.19}$$

Note that $A_1^q(\xi; m, n)$ is only non-zero for $m + n = 1$, whilst $A_2^q(\xi; m, n)$ and $A_3^q(\xi; m, n)$ are only non-zero for $m + n = 0$. We now define

$$\begin{aligned}
B_1^q(\xi; m, n) := & -\sigma_{m-1,n} + 4A\omega_{m,n} - \tau_{m+1,n} \\
& + \tau_{q-m+1,-q-n} - 4A\omega_{q-m,-q-n} + \sigma_{q-m-1,-q-n}, \text{ and}
\end{aligned} \tag{4.20}$$

$$\begin{aligned}
B_2^q(\xi; m, n) := & \tau_{m,n+1} - 4B\omega_{m,n} + \sigma_{m,n-1} \\
& - \tau_{q-m,-q-n+1} + 4B\omega_{q-m,-q-n} - \sigma_{q-m,-q-n-1}.
\end{aligned} \tag{4.21}$$

**Lemma 4.11**

For all $q$, and for all finite $\xi$, $A_1^q(\xi; m, n) = 0$ for all $(m, n)$ if and only if, for all $(m, n)$,

$$B_1^q(\xi; m, n) = B_2^q(\xi; m, n) = 0. \tag{4.22}$$

**Proof:** Indeed, for all $(m, n)$,

$$A_1^q(\xi; m, n) = B_1^q(\xi; m, n - 1) + B_2^q(\xi; m - 1, n) \ ,$$

and it follows that if (4.22) holds for all $(m,n)$, then $A_1^q(\xi;m,n) = 0$ for all $(m,n)$. Conversely, for all $(m,n)$,

$$\sum_{k=0}^{\infty} A_1^q(\xi;m-2k,n+2k+1) = \sum_{k=0}^{\infty}\big(-4B\omega_{m-2k-1,n+2k+1} + 4A\omega_{m-2k,n+2k}$$

$$- 4A\omega_{q-m+2k,-q-n-2k} + 4B\omega_{q-m+2k+1,-q-n-2k-1}\big)$$

$$- \tau_{m+1,n} + \sigma_{q-m-1,-q-n} ,$$

where the sums on the left- and right-hand side are in fact finite. Likewise, for all $(m,n)$,

$$\sum_{k=0}^{\infty} A_1^q(\xi;q-m+2k+2,-q-n-2k-1)$$

$$= \sum_{k=0}^{\infty}\big(-4B\omega_{q-m+2k+1,-q-n-2k-1} + 4A\omega_{q-m+2k+2,-q-n-2k-2}$$

$$- 4A\omega_{m-2k-2,n+2k+2} + 4B\omega_{m-2k-1,n+2k+1}\big)$$

$$+\tau_{q-m+1,-q-n} - \sigma_{m-1,n} .$$

Combining these identities yields

$$\sum_{k=0}^{\infty} A_1^q(\xi;m-2k,n+2k+1) + \sum_{k=0}^{\infty} A_1^q(\xi;q-m+2k+2,-q-n-2k-1) = B_1^q(\xi;m,n) ,$$

so that if $A_1^q(\xi;m,n)$ vanishes for all $(m,n)$, then so too does $B_1^q(\xi;m,n)$ for all $(m,n)$. In a similar manner, we show that if $A_1^q(\xi;m,n)$ vanishes for all $(m,n)$, then so too does $B_2^q(\xi;m,n)$ for all $(m,n)$, and this completes the proof. □

**Lemma 4.12**

For each $i \in \{1,2\}$, if $B_i^q(\xi;m,n) = 0$ for all $(m,n)$, then, for all $k \geqslant 0$,

$$B_i^{0,1}(\Theta_k^q(\xi);m,n) = 0 . \tag{4.23}$$

**Proof:** It suffices to address the case where $i = 1$, as the case where $i = 2$ is identical. Note first that $B_1^q(\xi;m,n)$ is only non-zero for $m+n = 0$. Next, for all $m < 0$,

$$B_1^{0,1}(\Theta_k^q(\xi);m,-m) = B_1^q(\xi;(m+k),-(m+k)) = 0 .$$

We likewise show that $B_1^{0,1}(\Theta_k^q(\xi);m,-m) = 0$ for all $m > 0$. We verify by inspection that

$$B_1^{0,1}(\mu;0,0) = 0$$

for all $\mu$, so that, in particular, $B_1^{0,1}(\Theta_k^q(\xi);0,0) = 0$, and the result follows. □

**Lemma 4.13**

For each $i \in \{2,3\}$, if $A_i^q(\xi;m,n) = 0$ for all $(m,n)$, then, for all $k \geqslant 0$, and for all $(m,n)$,

$$A_i^{0,1}(\Theta_k^q(\xi);m,n) = 0 \ . \tag{4.24}$$

**Proof:** It suffices to address the case where $i = 2$, as the case where $i = 3$ is identical. Note first that $A_2^q(\xi;m,n)$ is only non-zero for $m+n = 0$. Next, for $m < -1$,

$$A_2^{0,1}(\Theta_k^q(\xi);m,-m) = A_2^q(\xi;(m+k),-(m+k)) = 0 \ .$$

Likewise, we show that, for all $m > -1$,

$$A_2^{0,1}(\Theta_k^q(\xi);m,-m) = 0 \ .$$

We verify by inspection that, for all $\mu$,

$$A_2^{0,1}(\mu;0,0) = 0 \ ,$$

so that, in particular,

$$A_2^{0,1}(\Theta_k^q(\xi);0,0) = 0 \ .$$

It thus only remains to address the cases where $(m,n) \in \{(-1,1),(1,-1)\}$. However, denoting $\xi' := \Theta_k^q(\xi)$,

$$A_2^{0,1}(\Theta_k^q(\xi);-1,1) = -\omega_{k-2,2-k} - 4B\sigma_{k-2,1-k} + 4A\sigma_{k-1,-k} + \frac{1}{2}\big(\omega_{k,-k} + \omega_{q-k,-q+k}\big)$$
$$+ \frac{1}{2}\big(\omega_{k,-k} + \omega_{q-k,-q+k}\big) + 4\sigma_{q-k,-q+k-1}$$
$$- 4A\sigma_{q-k+1,-q+k-2} - \omega_{q-k+2,-q+k-2}$$
$$= A_2^q(\xi;k-1,1-k) \ ,$$

so that

$$A_2^{0,1}(\Theta_k^q(\xi);-1,1) = 0 \ .$$

In a similar manner, we show that

$$A_2^{0,1}(\Theta_k^q(\xi);1,-1) = 0 \ ,$$

and this completes the proof. $\square$

We now complete the proof of Theorem 4.9.

**Proof of Theorem 4.9:** Let $\text{Pot}(\Sigma)$ denote the manifold of potentials of $\Sigma$, and let $\text{Pot}_{\text{Skl}}(\Sigma)$ denote the submanifold consisting of those potentials which satisfy the $(g-1)$-Sklyanin condition. Define $\tilde{A} : \text{Pot}(\Sigma) \to \text{Jac}(\Sigma)$ such that, for all $\xi$, $\tilde{A}(\xi)$ is the image under the Abel map of the divisor of $\xi$. Note that $\tilde{A}$ is a smooth diffeomorphism onto its image. Let $\xi_0$ be a point of $\text{Pot}_{\text{Skl}}(\Sigma)$, denote $z_0 := \tilde{A}(\xi_0)$, and let $\text{Skl}(z_0)$ denote the Sklyanin subspace containing $z_0$. It will suffice to show that, for all such $\xi_0$, $D\tilde{A}(\xi_0)^{-1}$ maps the tangent space of $\text{Skl}(z_0)$ at $z_0$ isomorphically into the tangent space of $\text{Pot}_{\text{Skl}}(\Sigma)$ at $\xi_0$. Indeed, it then follows by the inverse function theorem that if $z_0 \in \text{Jac}(\Sigma)$ lies in the image under $\tilde{A}$ of $\text{Pot}_{\text{Skl}}(\Sigma)$, then so too does a neighbourhood of $z_0$ in $\text{Skl}(z_0)$, and a connectedness argument then allows us to conclude that the same holds for every non-special point of $\text{Skl}(z_0)$ with divisor supported in $\Sigma^*$. Furthermore, since every point of the real locus is non-special, the final assertion will also follow immediately.

We first show that every tangent vector to $\text{Skl}(z_0)$ at $z_0$ has the form $V[P,\overline{P}]$ for some complex polynomial $P$. Indeed, let $u$ be a tangent vector to $\text{Skl}(z_0)$ at $z_0$. By Lemma 2.7 and the subsequent remark, there exists a complex polynomial $P$ of order at most $(g-1)$ such that

$$u = V[2P,0] \ .$$

This vector is tangent to $\mathrm{Skl}(z_0)$ if and only if, for all $i$,

$$\mathrm{V}[P, 0]_i - \mathrm{V}[P, 0]_{(g+1)-i} = 0 \ .$$

By (2.15), this holds if and only if

$$\mathrm{V}[P, 0] = \mathrm{V}[0, \overline{P}] \ ,$$

so that

$$u = \mathrm{V}[P, 0] + \mathrm{V}[P, 0] = \mathrm{V}[P, 0] + \mathrm{V}[0, \overline{P}] = \mathrm{V}[P, \overline{P}] \ .$$

Conversely, we readily verify that every vector of this form is tangent to $\mathrm{Skl}(z_0)$, and the assertion follows.

Suppose now that

$$u := \mathrm{V}[e^{i\theta}\lambda^k, e^{-i\theta}\lambda^k] \ ,$$

for some real $\theta$, and some postive, integer $k$. By Lemma 3.2,

$$D\tilde{A}(\xi)^{-1} \cdot u = [\hat{\xi}(e^{i\theta}\lambda^k, e^{-i\theta}\lambda^k), \xi] \ .$$

However,

$$\hat{\xi}[e^{i\theta}\lambda^k, e^{-i\theta}\lambda^k] = e^{i\theta}\Pi_-(\lambda^{-k}\xi) + e^{i\theta}\Pi_+(\lambda^{k-(g-1)}\xi) = e^{i\theta}\Theta_k^{g-1}(\xi) \ .$$

By Lemma 4.10, $\Theta_k^{g-1}(\xi)$ satisfies the $(0,1)$-Sklyanin condition, from which it follows that $D\tilde{A}(\xi_0)^{-1} \cdot u$ satisfies the $(g-1)$-Sklyanin condition, and is therefore tangent to $\mathrm{Pot}_{\mathrm{Skl}}(\Sigma)$, and this completes the proof. $\square$

**4.3 - Proofs of main results.** We first show that Sklyanin subspaces in fact arise in complementary pairs.

**Lemma 4.14**

If (4.3) is satisfied over the Sklyanin subspace $\mathrm{Skl}(S_0)$, then the complementary condition

$$\phi(\lambda, \nu) = \frac{(4A - 4B\lambda)}{(\lambda - 1/\lambda)} \ . \tag{4.25}$$

is satisfied over its complementary Sklyanin subspace $\mathrm{Skl}(S \setminus S_0) = -\mathrm{Skl}(S_0)$.

In particular, if points of $\mathrm{Skl}(S_0)$ correspond to solutions of (1.1) satisfying the Durham boundary conditions (1.2), then points of $\mathrm{Skl}(S \setminus S_0) = -\mathrm{Skl}(S_0)$ correspond to solutions of (1.1) satisfying the complementary Durham boundary conditions (1.10).

**Proof:** Let $U$ denote the set of potential points of $\mathrm{Skl}(S_0)$. By hypothesis, $U$ is open and dense in $\mathrm{Skl}(S_0)$. Choose $z \in U$, let $D_z$ denote its preimage under the Abel map, let $\xi_z$ denotes its potential, and let $\phi_z$ denote its projectivized eigenline. By hypothesis, $\phi_z$ satisfies (4.3) at every point of $S_0$. Consider now the gauge transformed potential

$$\xi'_z := \left[ \begin{pmatrix} 0 & 1 \\ 1 & 0 \end{pmatrix} \xi \begin{pmatrix} 0 & 1 \\ 1 & 0 \end{pmatrix} \right]^t \ .$$

We verify that its projectivized eigenline is $\phi'_z := -(\phi_z \circ \sigma_1)$, its pole divisor is $D'_z := \sigma_1(D_z)$, and the image of $D'_z$ under the Abel map is $z' := -z$, so that

$$\Psi(\mathcal{A}([\phi'_z]_{\infty, fin})) = -\mathcal{A}(S_0) \ .$$

Since $S - (0)^8$ is the pole divisor of the meromorphic function $\lambda^{-2}\mathrm{Det}(K(\lambda))$, and since, by definition, $\mathcal{A}(0) = 0$,

$$\mathcal{A}(S) = 8\mathcal{A}((0)) = 0 \ ,$$

so that

$$\Psi(\mathcal{A}([\phi'_z]_{\infty, fin})) = \mathcal{A}(S) - \mathcal{A}(S_0) = \mathcal{A}(S \setminus S_0) \ .$$

It follows that $z'$ lies in the complementary Sklyanin subspace $\mathrm{Skl}(S \setminus S_0)$. Furthermore, $z'$ is the image under the Abel map of a non-special divisor supported on $\Sigma \setminus \{0, \infty\}$, and is therefore a potential point. Finally, we readily verify that $\phi'_z$ satisfies (4.25), and the result now follows by Theorem 4.9. $\square$

We are now ready to prove Theorems 1.1 and 1.2. For every special Sklyanin subspace $S_0$, we denote
$$\mathrm{Skl}_{\mathrm{re}}(S_0) := \mathrm{Skl}(S_0) \cap \mathrm{Jac}_{\mathrm{re}}(\Sigma) \ ,$$
and we call this subspace the *real special Sklyanin subspace* of $S_0$. Note, in particular, that since $\mathrm{Skl}(S_0)$ is a complex $g$-dimensional affine subspace of $\mathrm{Jac}(\Sigma)$, and since $\mathrm{Jac}_{\mathrm{re}}(\Sigma)$ is a real affine subspace of maximal dimension, $\mathrm{Skl}_{\mathrm{re}}(S_0)$ is a real $g$-dimensional subspace of $\mathrm{Jac}_{\mathrm{re}}(\Sigma)$.

**Proof of Theorem 1.1:** By Theorem 4.5, the Durham locus is contained within the union of all Sklyanin subspaces. By Theorem 4.9, the Durham locus is a union of Sklyanin subspaces. In particular, the real Durham locus is a union of real Sklyanin subspaces. It remains to show that any two connected components of the Durham locus are pairwise non-complementary. However, if one real Sklyanin subspace lies on the real Durham locus of the boundary condition (1.2), then, by Lemma 4.14, the complementary real Sklyanin subspace is contained within the real Durham locus of the complementary boundary condition (1.10). From this it follows that any two connected components of the real Durham locus are pairwise non-complementary, as desired. In particular, the real Durham locus may contain up to 2 connected components, and this completes the proof. $\square$

**Proof of Theorem 1.2:** Indeed, $z_0 + \mathrm{V}[iL, -iL]$ is an element of $\mathrm{Skl}^*_{\mathrm{re}}(\Sigma)$ if and only if
$$\Psi(z_0 + \mathrm{V}[iL, -iL]) = \mathcal{A}(S \setminus S_0) = -\mathcal{A}(S_0) \ .$$
However, since $z_0 \in \mathrm{Dur}_{\mathrm{re}}(\Sigma)$,
$$\mathcal{A}(S_0) = \Psi(z_0) \ .$$
Upon combining these relations, we see that $z_0 + \mathrm{V}[iL, -iL]$ is an element of $\mathrm{Dur}^*_{\mathrm{re}}(\Sigma)$ if and only if
$$2\Psi(z_0) + \Psi(\mathrm{V}[iL, -iL]) = 0 \ ,$$
as desired. $\square$

## A - Surfaces of constant boundary angle.

Let $u : \mathbb{R} \times [0, L] \to \mathbb{R}^3$ be a conformal immersion of constant mean curvature equal to $1/2$, where here we take the mean curvature to be equal to the algebraic mean of the principal curvatures. Let $\omega$ denote the conformal factor of its induced metric, let $\phi$ denote its Hopf differential, and suppose that
$$\phi = \frac{1}{4} dz dz \ .$$
Recall that the second fundamental form of $u$ is given by
$$\mathrm{II} = \phi + \frac{1}{2} e^{2\omega}(dx^2 + dy^2) + \overline{\phi} \ . \tag{A.1}$$
Let $\kappa_x$ and $\kappa_y$ denote the respective principal curvatures in the $x$ and $y$ directions. By (A.1),
$$\kappa_x = \frac{1}{2}(1 + e^{-2\omega}) \text{ and } \kappa_y = \frac{1}{2}(1 - e^{-2\omega}) \ . \tag{A.2}$$
Let $c_x : \mathbb{R} \times I \to \mathbb{R}$ denote the geodesic curvature of the horizontal lines with respect to the upward-pointing unit normal vector field.

**Lemma A.1**

*The geodesic curvature of horizontal lines satisfies*
$$c_x = e^{-\omega} \omega_y \ . \tag{A.3}$$

**Proof:** Indeed, the length element of horizontal lines is
$$dl_x = e^\omega dx \ .$$
The upward-pointing unit normal vector field is
$$\nu_x := e^{-\omega} \partial_y \ .$$
The geodesic curvature thus satisfies
$$c_x dl_x = D_{\nu_x} dl_x = e^{-\omega} D_{\partial_y} dl_x = e^{-\omega} \omega_y dl_x \ ,$$
and the result follows. $\square$

Consider now the restriction of $u$ to the horizontal line $X := \{0\} \times \mathbb{R}$.

**Theorem A.2**

*The restriction $u|_X$ lies on the surface of a sphere of radius $R$, and $u$ makes a constant angle $\theta$ with this sphere along this curve, if and only if, over $X$,*

$$\omega_y = \left(\frac{\varepsilon}{\sin(\theta)R} + \frac{\cos(\theta)}{2\sin(\theta)}\right)e^\omega + \frac{\cos(\theta)}{2\sin(\theta)}e^{-\omega} \ , \tag{A.4}$$

*for some $\varepsilon \in \{\pm 1\}$.*

**Proof:** Suppose first that $u|_X$ lies on the surface of a sphere of radius $R$ and that $u$ makes a constant angle $\theta$ with this sphere along this curve. Let $N : \mathbb{R} \times [0,L] \to \mathbb{S}^2$ denote the unit normal vector field over $u$. By hypothesis, at every point of $X$, the vector field

$$\xi := \sin(\theta)e^{-\omega}\partial_y u - \cos(\theta)N$$

is the unit normal vector field of a sphere of radius $R$. In particular, since $u(X)$ lies along this sphere, at every point of $X$,

$$\langle D_{\partial_x}\xi, \partial_x u\rangle = \frac{\varepsilon}{R}e^{2\omega} \ ,$$

for some $\varepsilon \in \{\pm 1\}$. However, differentiating the formula for $\xi$ in the $x$-direction yields, along $X$,

$$D_{\partial_x}\xi = -\sin(\theta)e^{-\omega}\omega_x\partial_y u + \sin(\theta)e^{-\omega}D_{\partial_x}\partial_y u - \cos(\theta)D_{\partial_x}N \ .$$

By conformality,

$$\langle \partial_y u, \partial_x u\rangle = 0 \ ,$$

by definition of the geodesic curvature,

$$\langle D_{\partial_x}\partial_y u, \partial_x u\rangle = e^{3\omega}c_x = e^{2\omega}\omega_y \ ,$$

and by definition of the principal curvature,

$$\langle D_{\partial_x}N, \partial_x u\rangle = e^{2\omega}\kappa_x = \frac{1}{2}e^{2\omega}(1+e^{-2\omega}) \ .$$

Combining these relations yields

$$\frac{\epsilon}{R}e^{2\omega} = e^\omega\sin(\theta)\omega_y - \frac{1}{2}e^{2\omega}\cos(\theta)(1+e^{-2\omega}) \ ,$$

and (A.4) follows.

Conversely, suppose that $\omega$ satisfies (A.4) over $X$. It suffices to show that the function

$$Y(x) := u(x) - R\sin(\theta)e^{-\omega}\partial_y u + R\cos(\theta)\nu$$

is constant over this line. Note, however, that for all $x$, the triplet $(\nu(x), e^{-\omega(x)}(\partial_x u)(x), e^{-\omega(x)}(\partial_y u)(x))$ is an orthonormal basis of $\mathbb{R}^3$. We verify by inspection that each of the components of $(D_{\partial_x}Y)(x)$ with respect to this basis vanishes. It follows that $D_{\partial_x}Y$ vanishes, so that $Y$ is constant, and this completes the proof. $\square$

## B - Bibliography.

[1] Andrews B., Li H., *Embedded constant mean curvature tori in the three-sphere*, J. Diff. Geom., **99**, no. 2, (2015), 169–189

[2] Belokolos E. D., Bobenko A. I., Enol'skii V. Z., Its A. R., Matveev V. B., *Algebro-Geometric Approach to Nonlinear Integrable Equations*, Springer Series in Nonlinear Dynamics, Springer-Verlag, (1994)

[3] Bobenko A. I., *All constant mean curvature tori in $\mathbb{R}^3$, $\mathbb{S}^2$, $\mathbb{H}^3$ in terms of theta functions*, Math. Ann., **290**, no. 2, (1991), 209–245

[4] Bobenko A., Kuksin S. B., *Small-amplitude solutions of the sine-Gordon equation on an interval under Dirichlet or Neumann boundary conditions*, J. Nonlinear. Sci., **5**, no. 3, (1995), 207–232

[5] Brendle S., *Embedded minimal tori in $\mathbb{S}^3$ and the Lawson conjecture*, Acta. Math., **211**, no. 2, (2013), 177–190

[6] Burstall F. E., Pedit F., *Harmonic maps via Adler–Kostant-Symes theory*, in *Harmonic Maps and Integrable Systems*, Fordy A. P., Wood J. C., (eds.), (Aspects of Mathematics), **23**, (1994)

[7] Carberry E., *Minimal tori in $\mathbb{S}^3$*, Pac. J. Math., **233**, no. 1, (2007), 41–69

[8] Cerezo A., Fernández I., Mira P., *Annular solutions to the partitioning problem in a ball*,

[9] Corrigan E., Dorey P. E., Rietdijk R. H., Sasaki R., Affine Toda field theory on a half line, *Phys. Lett. B*, **333**, (1994)

[10] Corrigan E., Dorey P. E., Rietdijk R. H., Aspects of addine Toda field theory on a half line, *Prog. Theor. Phys. Suppl.*, **118**, (1995), 143–164

[11] Delius G. W., Soliton-preserving boundary condition in affine Toda field theories, *Physics Letters B*, **444**, nos. 1–2, (1998), 217–223

[12] Ercolani N. M., Knörrer H., Trubowitz E., Hyperelliptic curves that generate constant mean curvature tori in R3, in *Integrable Systems (Luminy, 1991)*, Babelon O. et al. (eds.), Progress in Mathematics, **115**, Birkhaüser, Boston, 81–114, (1993)

[13] Farkas H. M., Kra I., *Riemann Surfaces*, Graduate Texts in Mathematics, **71**, Springer-Verlag, (1992)

[14] Fraser A., Schoen R., The first Steklov eigenvalue, conformal geometry, and minimal surface, *Adv. Math.*, **226**, no. 5, (2011), 4011–4030

[15] Fraser A., Schoen R., Sharp eigenvalue bounds and minimal surfaces in the ball, *Inv. Math.*, **203**, no. 3, (2016), 823–890

[16] Hauswirth L., Kilian M., Schmidt M. U., Finite type minimal annuli in $\mathbb{S}^2 \times \mathbb{R}$, *Illinois J. Math.*, **57**, no. 3, (2013), 697–741

[17] Jaggy C., On the classification of constant mean curvature tori in $\mathbb{R}^3$, *Comment. Math. Helv.*, **69**, no. 4, (1994), 640–658

[18] Kilian M., Schmidt M. U., On the moduli of constant mean curvature cylinders of finite type in the 3-sphere, arXiv:0712.0108

[19] Kilian M., Smith G., On the elliptic sinh-Gordon equation with integrable boundary conditions, *Nonlinearity*, **34**, no. 8, (2021), 5119–5135

[20] Klein S., *A spectral theory for simply-periodic solutions of the sinh-Gordon equation*, Lecture Notes in Mathematics, **2229**, Springer, (2018)

[21] Miranda R., *Algebraic Curves and Riemann Surfaces*, Graduate Studies in Mathematics, **5**, AMS, (1995)

[22] Pinkall U., Sterling I, On the Classification of Constant Mean Curvature Tori, *Ann. of Math.*, **130**, no. 2, (1989), 407-451

[23] Sklyanin E. K., Boundary conditions for integrable equations, *Funct. Anal. Appl.*, **21**, (1987), 164–166

[24] Sklyanin E. K., Boundary conditions for integrable quantum systems, *J. Phys. A. Math. Gen.*, **21**, (1988), 2375–2389

[25] Wente H. C., Counterexample to a conjecture of H. Hopf, *Pacific J. Math.*, **121** (1986), 193–243

[26] Wente H. C., Tubular capillary surfaces in a convex body, in *Advances in geometric analysis and continuum mechanics*, Concus P. and Lancaster K. (eds.), International Press, 288–298, (1995)